\documentclass[12pt,draft]{article}
\usepackage{amssymb}
\usepackage[top=72pt,bottom=96pt,left=72pt,right=72pt]{geometry}
\def\#{\sharp}
\def\<#1,#2>{\langle\,#1,\,#2\,\rangle}
\def\a{\mathfrak{a}}
\def\ad{\mathrm{ad}}
\def\b{\flat}
\def\C{\mathbb{C}}
\def\det{\mathrm{det}}
\def\dim{\mathrm{dim}\,}
\def\End{\mathrm{End}\,}
\def\FS{\mathrm{FS}}
\def\g{\mathfrak{g}}
\def\Gr{\mathrm{Gr}}
\def\H{\mathbb{H}}

\def\id{\mathrm{id}}

\def\k{\mathfrak{k}}
\def\Kr{\mathrm{Kr}}

\def\N{\mathbb{N}}
\def\O{\mathbb{O}}
\def\p{\mathfrak{p}}
\def\pf{\mathrm{pf}}
\def\pr{\mathrm{pr}}
\def\R{\mathbb{R}}
\def\Ric{\mathrm{Ric}}
\def\sec{\mathrm{sec}}
\def\S{\mathrm{Sym}\,}

\def\span{\mathrm{span}}
\def\St{\mathrm{St}}
\def\tr{\mathrm{tr}}
\def\vol{\mathrm{vol}}
\def\Vol{\mathrm{Vol}}
\def\X{\mathfrak{X}}
\def\Y{\mathfrak{Y}}
\def\Z{\mathbb{Z}}
\def\pfill{\par\vskip10pt plus3pt minus3pt\noindent}
\def\proof{\pfill\textbf{Proof:}\quad}
\def\qed{\ensuremath{\hfill\Box}}
\newtheorem{Lemma}{Lemma}[section]
\newtheorem{Remark}[Lemma]{Remark}
\newtheorem{Theorem}[Lemma]{Theorem}

\newtheorem{Corollary}[Lemma]{Corollary}

\newenvironment{Reference}[1]{\pfill\textbf{#1} \textit\bgroup}{\egroup\par}
\begin{document}
\title{Moments of Sectional Curvature}
\author{Gregor Weingart\footnote{
 address:\ Instituto de Matem\'aticas, Universidad Nacional Aut\'onoma
 de M\'exico, Avenida Universidad s/n, Colonia Lomas de Chamilpa, 62210
 Cuernavaca, Morelos, MEXIQUE; email: \texttt{gw@matcuer.unam.mx}.}}
\maketitle
\begin{center}
 \textbf{Abstract}
 \\[11pt]
 \parbox{412pt}{%
  The sectional curvature of a compact Riemannian manifold $M$ can be seen
  as a random variable on the Grassmann bundle of $2$--planes in $TM$ endowed
  with the Fubini--Study volume density. In this article we calculate the
  moments of this random variable by integrating suitable local Riemannian
  invariants and discuss the distribution of the sectional curvature of
  Riemannian products. Moreover we calculate the moments and the distribution
  of the sectional curvature for all compact symmetric spaces of rank $1$
  explicitly and derive a formula for the moments of general symmetric spaces.
  Interpolating the explicit values for the moments obtained we prove a weak
  version of the Hitchin--Thorpe Inequality.}
 \\[11pt]
\end{center}
\section{Introduction}
 Analyzing the sectional curvature is an indispensable tool in understanding
 the local geometry and the global topology of Riemannian manifolds. Given
 the importance of the sectional curvature we propose in this article to study
 the sectional curvature of a compact Riemannian manifold $M$ as a random
 variable defined on the Grassmann bundle of $2$--planes in $TM$. Our main
 motivation for this study besides pure curiosity is to provide us with a
 large subalgebra of local Riemannian invariants, whose values are easily
 calculated for many interesting examples and can be interpreted directly
 in terms of the underlying geometry. The more explicit values of local
 Riemannian invariants are known the easier it is to find and prove general
 linear relations between these invariants.

 The principal characteristic of the sectional curvature of a compact
 Riemannian manifold $M$ thought of as a random variable $\sec:\,\Gr_2TM
 \longrightarrow\R,\,(\,p,\,\Sigma\,)\longmapsto\sec_{R_p}(\Sigma),$ on
 the Grassmann bundle $\Gr_2TM$ of $2$--planes in $TM$ is of course its
 distribution, a probability measure on the real line, whose support
 determines the minimum and maximum of the sectional curvature of $M$.
 Treating the sectional curvature as a random variable may thus be seen
 as an alternative route to the calculation of pinching constants compared
 for example to the work of P\"uttmann \cite{p}. Determining the distribution
 of the sectional curvature directly however seems exceedingly difficult,
 for this reason we will focus instead on the {\em moments} of the sectional
 curvature with respect to a suitable volume density $\vol_\FS$ on $\Gr_2TM$.
 In essence the moments of the sectional curvature are integrated local
 Riemannian invariants:

 \begin{Reference}{Theorem \ref{gins}\ (Moments of Sectional Curvature)}
 \hfill\break
  Consider the category of euclidean vector spaces of dimension $m$ under
  isometries as morphisms and the functor, which associates to every such
  vector space $V$ with scalar product $g$ the vector space of algebraic
  curvature tensors $R:\,V\times V\times V\longrightarrow V,\,(X,Y,Z)
  \longmapsto R_{X,Y}Z$. The sequence $\{\,\Psi_k\,\}_{k\,\in\,\N_0}$
  of natural polynomials of degree $k$ on algebraic curvature tensors
  $$
   \Psi^V_k(\;R\;)
   \;\;:=\;\;
   \left.\frac{(\,-\,\Delta_g\,)^k}{[\,m\,+\,2k\,-\,2\,]_{2k}}\,\right|_0
   \left(\,X\,\longmapsto\,\exp\Big(\,\sum_{r>0}\frac1{2r}\,
   \tr_V(\,R_{\,\cdot\,,\,X}X\,)^r\,\Big)\,\right)
  $$
  defined using the positive Laplace operator $\Delta_g$ on $C^\infty V$
  associated to the scalar product $g$ and the falling factorial polynomial
  $[\,x\,]_s\,:=\,x(x-1)\ldots(x-s+1)$ can be used to calculate the moments
  of the sectional curvature of every compact Riemannian manifold $M$ via:
  $$
   E_M[\,\sec^k\,]
   \;\;:=\;\;
   \frac1{\Vol\,\Gr_2TM}\int_{\Gr_2TM}\!\!\sec^k_{R_p}(\,\Sigma\,)
   \,\vol_\FS(\,p,\Sigma\,)
   \;\;=\;\;
   \frac1{\Vol\,M}\int_M\Psi^{T_pM}_k(\,R_p\,)\,\vol_g(\,p\,)
  $$
 \end{Reference}

 \pfill
 Using this presentation we calculate the moments of the sectional curvature
 for all the compact symmetric spaces of rank $1$ up to covering, besides the
 round spheres $S^m$ this family of examples comprises the complex and
 hyperbolic projective spaces $\C P^n$ and $\H P^n$ with $n\,\geq\,2$ and
 the exceptional Cayley projective plane $\O P^2$. The sequence of moments
 for these symmetric spaces all correspond to very simple probability
 measures on $\R$, namely
 \begin{eqnarray*}
  E_{\C P^n}[\;F(\,\sec\,)\;]
  &=&
  \hbox to 64pt{${\displaystyle\frac16\;\,{n-\frac32\choose n-2}}$\hfill}
  \int_1^4\hbox to50pt{$\displaystyle\sqrt{\frac{s-1}3}^{-1}$\hss}
  \hbox to64pt{$\displaystyle\Big(\frac{4-s}3\Big)^{n-2}$\hfil}F(\,s\,)\;ds
  \\[2pt]
  E_{\H P^n}[\;F(\,\sec\,)\;]
  &=&
  \hbox to 64pt{${\displaystyle\frac36\;{2n-\frac32\choose 2n-3}}$\hfill}
  \int_1^4\hbox to50pt{$\displaystyle\sqrt{\frac{s-1}3}$\hfil}
  \hbox to64pt{$\displaystyle\Big(\frac{4-s}3\Big)^{2n-3}$\hfil}F(\,s\,)\;ds
  \\[2pt]
  E_{\O P^2}[\;F(\,\sec\,)\;]
  &=&
  \hbox to 64pt{${\displaystyle\frac76\quad\;{\frac{13}2\choose 3}}$\hfill}
  \int_1^4\hbox to50pt{$\displaystyle\sqrt{\frac{s-1}3}^5$\hfil}
  \hbox to64pt{$\displaystyle\Big(\frac{4-s}3\Big)^3$\hfil}F(\,s\,)\;ds
 \end{eqnarray*}
 where in all cases the metric is normalized to have sectional curvatures
 in the interval $[\,1,\,4\,]$. For general compact symmetric spaces the
 calculation of the moments of the sectional curvature reduces via a suitable
 variant of Weyl's Integration Formula to an integration over spherical
 simplices, in this way we obtain closed formulas for the moments of the
 sectional curvature of the Fubini--Study metric on the Grassmannians of
 $2$--planes in $\R^n$, $\C^n$ and $\H^n$. In order to augment this limited
 stock of examples we consider the behaviour of the moments of the sectional
 curvature under Riemannian products and obtain from Theorem \ref{gins}:

 \begin{Reference}{Lemma \ref{mp}\ (Sectional Curvature Moments of Products)}
 \hfill\break
  Consider two compact Riemannian manifolds $M$ and $N$ of dimensions
  $m,\,n\,\geq\,0$. The moments of the sectional curvature of the Riemannian
  product $M\times N$ can be calculated via
  $$
   E_{M\times N}[\;\sec^k\;]
   \;\;=\;\;
   \sum_{r\,=\,0}^k{k\choose r}\,\frac{[\,m+2r-2\,]_{2r}
   \,[\,n+2(k-r)-2\,]_{2(k-r)}}{[\,m+n+2k-2\,]_{2k}}
   \,E_M[\;\sec^r\;]\,E_N[\;\sec^{k-r}\;]
  $$
  where $[\,x\,]_s\,:=\,x(x-1)\ldots(x-s+1)$ denotes the falling factorial
  polynomial as before.
 \end{Reference}

 \pfill
 From the point of view of combinatorics this product formula for the moments
 of the sectional curvature appears rather strange, nevertheless it has a very
 nice probabilistic model in form of a family of probability measures on the
 $2$--simplex $\Delta_2\,\subset\,\R^2$. This probabilistic model allows us
 to extend the product formula for moments to a product formula for the
 distribution of the sectional curvature of Riemannian products formulated
 in Corollary \ref{scx}.

 \pfill
 The preceding results indicate that the moments of the sectional curvature
 are best seen in the context of integrated local Riemannian invariants.
 Interestingly the algebra of all local Riemannian invariants allows for
 a graphical calculus similar to the graph algebras describing the
 Rozansky--Witten invariants \cite{rw} of hyperk\"ahler manifolds. According
 to Theorem \ref{gins} the moments of the sectional curvature generate a
 large and quite explicit subalgebra of local Riemannian invariants, which
 are easy to calculate and relate directly to Riemannian geometry. In a
 forthcoming publication \cite{w} we will combine our knowledge of explicit
 moments of sectional curvature together with the graphical calculus
 mentioned above in order to prove linear relations between local Riemannian
 invariants by interpolation. Interpolating for example the first and second
 moments of the sectional curvature of the Riemannian manifolds $\C P^2$,
 $S^4$, $S^2\times S^2$ and $S^1\times S^3$ we easily obtain the following
 version
 \begin{eqnarray*}
  \lefteqn{\frac{4\,\pi^2}{3}\;\chi(\,M\,)\;+\;
  \frac49\,\int_M|\;\Ric^\circ_p\;|^2\,\vol_g(\,p\,)}
  \quad
  &&
  \\
  &=&
  \int_M\Psi_2(\,R_p\,)\,\vol_g(\,p\,)\;+\;
  4\,\int_M\Big[\;\Psi_2(\,R_p\,)\;-\;\Psi_1(\,R_p\,)^2\;\Big]\,\vol_g(\,p\,)
 \end{eqnarray*}
 of the Hitchin--Thorpe Inequality \cite{th} for compact $4$--dimensional
 Riemannian manifolds $M$, where $\chi(M)$ is the Euler characteristic and
 $\Ric^\circ\,:=\,\Ric-\frac\kappa4\,\id$ the trace free Ricci tensor.

 \pfill
 Organisatorically this article is divided into three sections besides this
 introduction. In Section \ref{ing} we reformulate the integration over
 Grassmannians into a double integration over round spheres and use this
 reformulation together with an integration trick for integrating polynomials
 over round spheres to prove Theorem \ref{gins}. Explicit examples for the
 moments and the distribution of sectional curvature are given in Section
 \ref{ex}, while Section \ref{mdp} studies Riemannian products and discusses
 the interpolation proof of the Hitchin--Thorpe Inequality. The author would
 like to thank U.~Semmelmann, the students and the staff of the University
 of Stuttgart for their hospitality enjoyed during several research stays
 in Stuttgart.
\section{Integration over Grassmannians}
\label{ing}
 In general the moments of a random variable are defined as integrals
 over the probability space in question, the moments of the sectional
 curvature of an algebraic curvature tensor are thus integrals over the
 Grassmannian of $2$--planes in the underlying euclidean vector space.
 Integration is in itself nothing else but the evaluation of a very
 special linear functional on an infinite dimensional vector space,
 hence it is prudent to look out for integration tricks before embarking
 on an explicit integration in local coordinates or otherwise. In this
 section we discuss such an integration trick for the integration of
 polynomials over a Grassmannian of $2$--planes, which we use in turn
 to encode the moments of the sectional curvature of an algebraic curvature
 tensor in a simple generating formal power series.

 \pfill
 Recall that an algebraic curvature tensor on a euclidean vector space
 $V$ of dimension $m$ with scalar product $g\,\in\,\S^2V^*$ is a trilinear
 map $R:\,V\times V\times V\longrightarrow V,\,(X,Y,Z)\longmapsto R_{X,Y}Z,$
 which is alternating $R_{X,Y}Z\,=\,-R_{Y,X}Z$ in its first two arguments,
 defines skew--symmetric endomorphisms $Z\longmapsto R_{X,Y}Z$ of $V$ for
 all $X,\,Y\,\in\,V$ and satisfies the first Bianchi identity:
 $$
  R_{X,\,Y}Z\;+\;R_{Y,\,Z}X\;+\;R_{Z,\,X}Y
  \;\;=\;\;
  0
 $$
 Every such algebraic curvature tensor $R$ is completely determined by is
 associated sectional curvature, which is a well--defined function
 on the Grassmannian $\Gr_2V$ of $2$--planes in $V$:
 $$
  \sec_R:\quad\Gr_2V\;\longrightarrow\;\R,\qquad
  \span\,\{\,X,\,Y\,\}\;\longmapsto\;\frac{g(\;R_{Y,\,X}X,\;Y\;)}
  {g(\,X,\,X\,)\,g(\,Y,\,Y\,)\;-\;g(\,X,\,Y\,)^2}
 $$
 In the sequel we want to integrate the powers $\sec_R^k$ of the sectional
 curvature of an algebraic curvature tensor $R$ against the volume density
 of the Grassmannian $\Gr_2V$ endowed with the Fubini--Study metric $g^\FS$.
 The principal idea to make this integration feasible is to convert integrals
 over $\Gr_2V$ into double integrals over unit spheres, because integrals
 over unit spheres are relatively easy to calculate, compare for example
 Corollary \ref{sphint}. The part of the Cardan joint of this conversion
 is played by the Stiefel manifold of orthonormal $2$--frames:
 $$
  \St_2V
  \;\;:=\;\;
  \{\;(X,Y)\;\;|\;\;g(X,X)\,=\,1\,=\,g(Y,Y)
  \;\textrm{\ and\ }\;g(X,Y)\,=\,0\;\}
 $$
 Evidently $\St_2V\,\subset\,V\times V$ is a submanifold with tangent space
 equal to
 \begin{equation}\label{stts}
  T_{(X,Y)}\St_2V
  \;\;=\;\;
  \{\;(\dot X,\dot Y)\;\;|\;\;
  g(X,\dot X)\,=\,0\,=\,g(Y,\dot Y)
  \;\textrm{\ and\ }\;g(\dot X,Y)\,+\,g(X,\dot Y)\,=\,0\;\}
 \end{equation}
 in an orthonormal $2$--frame $(X,Y)\,\in\,\St_2V$. In order to have $\St_2V$
 acting as a Cardan joint we need to define a Riemannian metric $g^\FS$ on
 $\St_2V$ so that the canonical projections
 $$
  \begin{array}{lccccccc}
   \pr_\Gr:&\St_2V&\!\longrightarrow\!&\Gr_2V,&\qquad&
   (\,X,\,Y\,)&\!\longmapsto\!&\span\,\{\;X,\,Y\;\}
   \\[3pt]
   \pr_S\;\,:&\St_2V&\!\longrightarrow\!&
   S_V\hbox to0pt{\hskip5pt,\hss}&\qquad&(\,X,\,Y\,)&\!\longmapsto\!&X
  \end{array}
 $$
 to $\Gr_2V$ and the unit sphere $S_V\,\subset\,V$ are Riemannian submersions
 with respect to the Fubini--Study metric on $\Gr_2V$ and the round metric
 $g$ on $S_V$. In terms of the description (\ref{stts}) of the tangent space
 $T_{(X,Y)}\St_2V$ in $(X,Y)\,\in\,\St_2V$ the Riemannian metric of choice
 reads:
 $$
  g^\FS_{(X,Y)}\Big(\;(\dot X_1,\dot Y_1),\;(\dot X_2,\dot Y_2)\;\Big)
  \;\;:=\;\;
  g(\;\dot X_1,\;\dot X_2\;)\;+\;g(\;\dot Y_1,\;\dot Y_2\;)
  \;-\;g(\;Y,\;\dot X_1\;)\,g(\;Y,\;\dot X_2\;)
 $$
 The interesting and rather unexpected final term in this definition can be
 written in a number of equivalent ways due to the relation $g(Y,\dot X)\,=\,
 -g(X,\dot Y)$ valid for all tangent vectors $(\dot X,\dot Y)\,\in\,T_{(X,Y)}
 \St_2V$. Of course we would get a Riemannian metric on $\St_2V$ even without
 this final term, the Fubini--Study metric $g^\FS$ however has the very
 convenient property that the tangent space $T_{(X,Y)}\St_2V$ decomposes
 orthogonally into the direct sums
 \begin{eqnarray*}
  \lefteqn{T_{(X,Y)}\St_2V}
  &&
  \\
  &=&
  \R\,(-Y,X)\;\oplus\;\{\;(\dot X,\dot Y)\;|\;\;
  \dot X,\,\dot Y\;\textrm{\ orthogonal to both\ }\;X\;\textrm{\ and\ }\;Y\;\}
  \\
  &=&
  \{\;(\,0\,,\dot Y)\;|\;\;\dot Y\;\textrm{\ orthogonal to\ }
  \;X,\,Y\;\}\;\oplus\;\{\;(\dot X,-g(\dot X,Y)X)\;|\;\;\dot X
  \;\textrm{\ orthogonal to\ }\;X\;\}
 \end{eqnarray*}
 where $\R\,(-Y,X)$ is the vertical tangent space for the projection to the
 Grassmannian $\Gr_2V$, while $\{\,(0,\dot Y)\,|\,\dot Y\textrm{\ orthogonal
 to\ }X,Y\,\}$ is the vertical tangent space for the projection $\pr_S$ to
 the unit sphere $S_V$. In consequence the projection $\pr_S$ is a Riemannian
 submersion, because:
 \begin{eqnarray*}
  \lefteqn{g^\FS_{(X,Y)}\Big(\;(\,\dot X_1,\,-g(\dot X_1,Y)\,X\,),
   \;(\,\dot X_2,\,-g(\dot X_2,Y)\,X\,)\;\Big)}
  \qquad
  &&
  \\
  &=&
  g(\,\dot X_1,\,\dot X_2\,)
  \;+\;g(\,\dot X_1,\,Y\,)\,g(\,\dot X_2,\,Y\,)
  \;-\;g(\,Y,\,\dot X_1\,)\,g(\,Y,\,\dot X_2\,)
  \;\;=\;\;
  g(\,\dot X_1,\,\dot X_2\,)
 \end{eqnarray*}
 For precisely this reason we had to include the final term into the
 definition of the Fubini--Study metric $g^\FS$ on the Stiefel manifold
 $\St_2V$. In passing we observe that the fiber of the projection $\pr_S$
 to the unit sphere $S_V$ in $X\,\in\,S_V$ is the round unit sphere
 $S_{\{\,X\,\}^\perp}$ in the orthogonal complement of $X$ due to
 $g^\FS_{(X,Y)}(\,(0,\dot Y_1),\,(0,\dot Y_2)\,)\,=\,g(\,\dot Y_1,
 \,\dot Y_2\,)$.

 In order to verify that the projection $\pr_\Gr$ to the Grassmannian
 $\Gr_2V$ of $2$--planes is likewise a Riemannian submersion we recall
 that the Fubini--Study metric $g^\FS$ is defined by means of the embedding
 of $\Gr_2V$ into the vector space of symmetric endomorphisms of $V$
 $$
  \Gr_2V\;\longrightarrow\;\End_{\mathrm{sym}}V,\qquad
  \Sigma\;\longmapsto\;\pr_\Sigma
 $$
 and the scalar product $h(F,\tilde F)\,:=\,\frac12\,\tr_V(F\tilde F)$
 on $\End_{\mathrm{sym}}V$. In terms of the musical isomorphism $\#:\,
 V\longrightarrow V^*,\,X\longmapsto g(X,\,\cdot\,),$ we may write
 $\pr_\Sigma\,=\,X^\#\otimes X\,+\,Y^\#\otimes Y$ for some orthonormal
 $2$--frame $(X,Y)\,\in\,\St_2V$ representing $\Sigma\,=\,\span\{\,X,\,
 Y\,\}\,\in\,\Gr_2V$ and find
 $$
  \dot\pr_\Sigma
  \;\;=\;\;
  \dot X^\#\otimes X\;+\;X^\#\otimes\dot X
  \;+\;\dot Y^\#\otimes Y\;+\;Y^\#\otimes\dot Y 
 $$
 for a tangent vector $(\dot X,\dot Y)\,\in\,T_{(X,Y)}\St_2V$. Calculating
 the trace scalar product results in:
 $$
  \frac12\,\tr_V(\,\dot\pr^2_\Sigma\,)
  \;\;=\;\;
  g(\dot X,\dot X)\,+\,g(\dot Y,\dot Y)
  \,+\,2\,g(Y,\dot X)\,g(X,\dot Y)
  \;\;=\;\;
  g(\pr^\perp_\Sigma\dot X,\dot X)\,+\,g(\pr^\perp_\Sigma\dot Y,\dot Y)
 $$
 In consequence we conclude that the projection $\pr_\Gr$ from the Stiefel
 manifold $\St_2V$ of orthonormal $2$--frames to the Grassmannian $\Gr_2V$
 of $2$--planes is a Riemannian submersion as well, its fibers are two
 disjoint circles of length $2\pi$ due to $g^\FS_{(X,Y)}((-Y,X),(-Y,X))\,=\,1$.

 The idea of using the Riemannian submersions $\pr_\Gr$ and $\pr_S$ to
 convert integrals over the Grassmannian $\Gr_2V$ into double integrals
 over unit spheres is illustrated by the calculation
 $$
  \Vol\;\St_2V
  \;\;=\;\;
  \frac{2\,\pi^{\frac m2}}{\Gamma(\,\frac m2\,)}\,
  \frac{2\,\pi^{\frac {m-1}2}}{\Gamma(\,\frac {m-1}2\,)}
  \;\;=\;\;
  \frac{2^m\,\pi^{m-1}}{(\,m-2\,)!}
  \qquad\qquad
  \Vol\;\Gr_2V
  \;\;=\;\;
  \frac{(\,2\,\pi\,)^{m-2}}{(\,m-2\,)!}
 $$
 of the volumes of the Stiefel manifold and the Grassmannian with respect
 to their Fubini--Study metrics, in which we use Euler's formula $\Vol\,S_V
 \,=\,\frac{2\,\pi^{\frac m2}}{\Gamma(\,\frac m2\,)}$ for the volume of unit
 spheres and simply divide by the volume $4\pi$ of two disjoint circles to
 obtain $\Vol\,\Gr_2V$.

 \pfill
 Turning away from the Stiefel manifold and the Grassmannian we want to
 discuss a couple of integration tricks, which make the calculation of
 integrals over unit spheres feasible. As before we consider a euclidean
 vector space $V$ of dimension $m$ with scalar product $g\,\in\,\S^2V^*$.
 The dual scalar product $g^{-1}$ on $V^*$ is defined by declaring the
 mutually inverse musical isomorphisms $\#:\,V\longrightarrow V^*,\,
 X\longmapsto g(X,\,\cdot\,),$ and $\b:\,V^*\longrightarrow V$ to be
 isometries, this is
 \begin{equation}\label{ginv}
  g^{-1}(\;\alpha,\;\beta\;)
  \;\;:=\;\;
  g(\;\alpha^\b,\;\beta^\b\;)
  \;\;=\;\;
  \alpha(\,\beta^\b\,)
  \;\;=\;\;
  \sum_{\mu\,=\,1}^m\alpha(\,X_\mu\,)\,\beta(\,X_\mu\,)
 \end{equation}
 where $X_1,\,\ldots,\,X_m$ is some orthonormal basis for $V$. An orthonormal
 basis $x_1,\,\ldots,\,x_m\,\in\,V^*$ of $V^*$ with respect to the scalar
 product $g^{-1}$ can be thought of as a global orthonormal coordinate system
 for $V$ considered as a manifold endowed with the translation invariant
 Riemannian metric $g$, in particular the Laplace--Beltrami operator on
 $C^\infty V$ becomes:
 \begin{equation}\label{lb}
  \Delta_g
  \;\;:=\;\;
  -\,\sum_{\mu\,=\,1}^m\frac{\partial^2}{\partial x_\mu^2}
 \end{equation}

 \begin{Lemma}[Integration against Gaussian Probability Density]
 \hfill\label{hk}\break
  Let $V$ be a euclidean vector space of dimension $m$ with scalar product
  $g\in\S^2V^*$ considered as a translation invariant Riemannian metric on
  $V$. The Laplace--Beltrami operator $\Delta_g$ associated to $g$ allows us
  integrate $p\,\in\,\S\,V^*$ against the Gaussian probability density via:
  $$
   \left.e^{-\frac12\,\Delta_g}\right|_0\;p
   \;\;=\;\;
   \frac1{(2\pi)^{\frac m2}}\;\int_V
   e^{-\,\frac12\,g(X,X)}\,p(\,X\,)\,\vol_g(\,X\,)
  $$
 \end{Lemma}

 \proof
 Consider the exponential function $X\longmapsto e^{\alpha(X)}$ associated
 to a given $\alpha\,\in\,V^*$. Completing the square this function can be
 integrated easily against the Gaussian probability density associated to the
 translation invariant Riemannian metric $g$ on $V$
 \begin{eqnarray*}
  \frac1{(2\pi)^{\frac m2}}\,\int_V
   e^{-\,\frac12\,g(X,X)}\,e^{\alpha(X)}\,\vol_g(\,X\,)
  &=&
  \frac1{(2\pi)^{\frac m2}}\,\int_V
   e^{-\,\frac12\,g(X-\alpha^\b,X-\alpha^\b)}
   \,e^{\frac12\,g^{-1}(\alpha,\alpha)}\,\vol_g(\,X\,)
  \\
  &=&
  e^{\frac12\,g^{-1}(\alpha,\alpha)}
  \;\;=\;\;
  \left.e^{-\frac12\,\Delta_g}\right|_0\;e^\alpha
 \end{eqnarray*}
 using the consequence $(-\Delta_g)\,e^\alpha\,=\,g^{-1}(\alpha,\alpha)\,
 e^\alpha$ of equation (\ref{ginv}). The stated integration formula thus
 holds true for the homogeneous pieces $p\,=\,\frac1{k!}\alpha^k$ of the
 generating function $e^\alpha$ of degrees $k\,\geq\,0$. By polarization
 and linearity this results extends to all of $\S\,V^*$.
 \qed

 \pfill
 In passing we remark that Lemma \ref{hk} reflects the characteristic property
 of the heat kernel of the Laplace--Beltrami operator on a flat Riemannian
 manifold, compare \cite{bgv}. In the calculation of integrals over the
 Stiefel manifold and in turn over the Grassmannian we will make use of
 two corollaries of Lemma \ref{hk}. Before discussing these corollaries
 however let us recall the definition of the falling factorial polynomial
 $[\,x\,]_s$ with parameter $s\,\in\,\Z$, namely $[\,x\,]_s
 \,:=\,x(x-1)\ldots(x-s+1)$ for positive $s\,>\,0$, while $[\,x\,]_0\,:=\,1$
 and $[\,x\,]_s\,:=\,0$ for $s\,<\,0$. The most important property of the
 falling factorial polynomial is the twisted symmetry $[\,x\,+\,s\,]_s\,=\,
 (-1)^s[\,-\,x\,-\,1]_s$ valid for all $s\,\in\,\Z$, moreover the generalized
 binomial coefficients ${x\choose s}\,:=\,\frac1{s!}\,[\,x\,]_s$ appear in
 Newton's power series expansion $(1+t)^x\,=\,\sum_{s\geq0}{x\choose s}\,t^s$
 valid for $|\,t\,|\,<\,1$ and all $x\,\in\,\R$. In turn they appear in the
 proof of the following corollary:

 \begin{Corollary}[Integration of Symmetric Endomorphisms]
 \hfill\label{igauss}\break
  Consider a symmetric endomorphism $F:\,V\longrightarrow V$ of a euclidean
  vector space $V$ of dimension $m$ with scalar product $g\,\in\,\S^2V^*$.
  The exponential function $\,X\longmapsto\exp(\,t\,g(FX,X)\,)$ can be
  integrated against the Gaussian probability density $(2\pi)^{-\frac m2}
  \,e^{-\frac 12\,g}\,\vol_g$ on $V$ using:
  $$
   \frac1{(2\pi)^{\frac m2}}\;\int_V
   e^{-\frac 12\,g(X,X)}\,\exp\Big(\,t\,g(FX,X)\,\Big)\,\vol_g(\,X\,)
   \;\;=\;\;
   \exp\Big(\;\sum_{r>0}\frac{(2t)^r}{2r}\,\tr_VF^r\;\Big)
  $$
 \end{Corollary}

 \proof
 According to Sylvester's Theorem of Inertia we can choose an orthonormal
 basis $X_1,\,\ldots,\,X_m$ of $V$ consisting of eigenvectors for the
 symmetric endomorphism $F$ with eigenvalues $f_1,\,\ldots,\,f_m\,\in\,\R$
 respectively. In terms of the dual orthonormal basis $x_1,\,\ldots,\,x_m$
 the function $\,X\longmapsto g(FX,X)$ reads $\sum_{\mu=1}^m f_\mu x_\mu^2$,
 whereas the Laplace--Beltrami operator $\Delta_g$ is given by (\ref{lb})
 as before. The integration formula of Lemma \ref{hk} thus becomes the
 formula:
 \begin{eqnarray*}
  \lefteqn{\frac1{(2\pi)^{\frac m2}}\;\int_V
  e^{-\frac 12\,g(X,X)}\,\exp\Big(\,t\,g(FX,X)\,\Big)\,\vol_g(\,X\,)}
  \qquad
  \\
  &=&
  \left.e^{-\frac12\,\Delta_g}\right|_0
  \prod_{\mu\,=\,1}^m\exp(\;tf_\mu x_\mu^2\;)
  \;\;=\;\;
  \prod_{\mu\,=\,1}^m\Big(\;\left.
  e^{\frac12\,\frac{\partial^2}{\partial x_\mu^2}}\right|_0
  \exp(\;tf_\mu x_\mu^2\;)\;\Big)
 \end{eqnarray*}
 The factors on the right hand can be evaluated explicitly using a little
 bit of combinatorics
 \begin{eqnarray*}
  \left.e^{\frac12\frac{\partial^2}{\partial x^2}}\right|_0
  \exp(\,tfx^2\,)
  &=&
  \sum_{r\geq 0}\frac{(2\,tf)^r}{4^r\,r!^2}\,
  \left.\frac{\partial^{2r}}{\partial x^{2r}}\right|_0
  x^{2r}
  \;\;=\;\;
  \sum_{r\geq 0}{-\frac12\choose r}\,(\,-2\,tf\,)^r
  \\
  &=&
  (1-2\,tf)^{-\frac12}
  \;\;=\;\;
  \exp\Big(-\frac12\,\log(1-2\,tf)\Big)
  \;\;=\;\;
  \exp\Big(\,\sum_{r>0}\frac{(2t)^r}{2r}\,f^r\,\Big)
 \end{eqnarray*}
 besides the triviality $\frac1{4^r}{2r\choose r}\,=\,(-1)^r{-\frac12
 \choose r}$ and Newton's power series. With $\exp$ being multiplicative
 the integration formula thus follows from $\tr_VF^r\,=\,\sum_{\mu\,=\,1}^m
 f_\mu^r$.
 \qed

 \begin{Corollary}[Integration of Polynomials over Spheres]
 \hfill\label{sphint}\break
  The integral of a homogeneous polynomial $p\in\S^kV^*$ of degree
  $k\,\in\,\N_0$ over the unit sphere $S_V\,\subset\,V$ of a euclidean
  vector space $V$ of dimension $m$ with scalar product $g$ can be
  evaluated by using the Laplace--Beltrami operator $\Delta_g$ associated
  to the translation invariant Riemannian metric $g$ on $V$. More precisely
  this integration trick tells us for even $k\,=\,2\kappa$:
  $$
   \frac1{\Vol\;S_V}\int_{S_V}p(\,X\,)\,\vol_g(\,X\,)
   \;\;=\;\;
   \left.\frac{(\,-\,\Delta_g\,)^\kappa}
   {4^\kappa\,\kappa!\,[\,\frac m2\,+\,\kappa\,-\,1\,]_\kappa}
   \right|_0\;p
  $$
  The antisymmetry of $p$ under the antipodal map makes the left hand side
  vanish for odd $k$.
 \end{Corollary}

 \proof
 The Laplace--Beltrami operator $\Delta_g$ maps homogeneous polynomials
 to homogeneous polynomials of course so that $\left.(-\Delta)^\kappa\right|_0
 p\,=\,0$ for every homogeneous polynomial $p\,\in\,\S^kV^*$ of degree $k$
 unless $k\,=\,2\kappa$ is even. By converting the integral over $V$ to
 polar coordinates we may thus rewrite Lemma \ref{hk} for a homogeneous
 polynomial $p$ of even degree
 \begin{eqnarray*}
  \left.e^{-\,\frac12\,\Delta_g}\right|_0p
  \;\;=\;\;
  \left.\frac{(\,-\,\Delta_g\,)^\kappa}{2^\kappa\,\kappa!}\right|_0\;p
  &=&
  \frac1{(2\pi)^{\frac m2}}\;\int_V
  e^{-\,\frac12\,g(X,X)}\,p(\,X\,)\,\vol_g(\,X\,)
  \\
  &=&
  \frac1{(2\pi)^{\frac m2}}
  \;\int_0^\infty e^{-\frac12r^2}\,r^{m+k-1}\,dr
  \;\int_{S_V}p(\,X\,)\,\vol_g(\,X\,)
  \\
  &=&
  \frac{2^\kappa}{\Vol\;S_V}\,
  [\,{\textstyle\frac m2}\,+\,\kappa\,-\,1\,]_\kappa
  \,\int_{S_V}p(\,X\,)\,\vol_g(\,X\,)
 \end{eqnarray*}
 using the substitution $\rho\,=\,\frac12r^2,\,d\rho\,=\,r\,dr,$ and the
 definition of the $\Gamma$--function to evaluate
 $$
  \int_0^\infty e^{-\frac12r^2}\,r^{m+k-1}\,dr
  \;\;=\;\;
  \int_0^\infty e^{-\rho}\,(2\rho)^{\frac{m+k-2}2}\,d\rho
  \;\;=\;\;
  2^{\kappa+\frac m2\,-\,1}\,
  [\,{\textstyle\frac m2}\,+\,\kappa\,-\,1\,]_\kappa\,
  \Gamma(\,{\textstyle\frac m2}\,)
 $$
 the radial integral as well as $\Gamma(x+1)\,=\,x\,\Gamma(x)$ and Euler's
 formula $\Vol\,S_V\,=\,\frac{2\,\pi^{\frac m2}}{\Gamma(\,\frac m2\,)}$ for
 the volume of the unit sphere in dimension $m$.
 \qed

 \pfill
 In passing we want to point out a special consequence of the preceding
 proof, the integral of a homogeneous polynomial $p\,\in\,\S^kV$ of even
 degree $k\,=\,2\kappa$ against the volume density on $S_V$ agrees up to
 a conversion factor, which only depends on $\kappa$ and the dimension $m$,
 with the integral of $p$ against the Gaussian probability density
 $(2\pi)^{-\frac m2}\,e^{-\frac12g}\,\vol_g$ on $V$:
 \begin{equation}\label{conv}
  \frac1{(2\pi)^{\frac m2}}\int_V
  e^{-\,\frac12\,g(X,X)}\,p(\,X\,)\,\vol_g(\,X\,)
  \;\;=\;\;
  2^\kappa\,[\,{\textstyle\frac m2}\,+\,\kappa\,-\,1\,]_\kappa
  \,\Big(\,\frac1{\Vol\;S_V}\int_{S_V}p(\,X\,)\,\vol_g(\,X\,)\,\Big)
 \end{equation}
 This reformulation of Corollary \ref{sphint} comes in handily in the proof
 of the following theorem:

 \begin{Theorem}[Moments of Sectional Curvatures]
 \hfill\label{gins}\break
  Let $R:\,V\times V\times V\longrightarrow V$ be an algebraic curvature
  tensor on a euclidean vector space $V$ of dimension $m$ with scalar product
  $g$. The powers $\sec^k_R,\,k\,\in\,\N_0,$ of the sectional curvature
  function $\sec_R:\,\St_2V\longrightarrow\R,\,(X,Y)\longmapsto
  g(R_{Y,X}X,Y),$ integrate over $\St_2V$ to the expression:
  \begin{eqnarray*}
   \lefteqn{\frac1{\Vol\;\St_2V}\int_{\St_2V}g(\;R_{Y,\,X}X,\;Y\;)^k
   \,\vol_\FS(\,X,\,Y\,)}\qquad\qquad
   &&
   \\
   &=&
   \left.\frac{(\,-\,\Delta_g\,)^k}{[\,m\,+\,2k\,-\,2\,]_{2k}}\,\right|_0
   \left(\,X\,\longmapsto\,\exp\Big(\,\sum_{r>0}\frac1{2r}\,
   \tr_V(\,R_{\,\cdot\,,\,X}X\,)^r\,\Big)\,\right)
  \end{eqnarray*}
  In particular we obtain for the scalar curvature $\kappa(\,R\,)$ of the
  algebraic curvature tensor $R$:
  $$
   \frac{\kappa(\,R\,)}{m\,(\,m\,-\,1\,)}
   \;\;=\;\;
   \frac1{\Vol\;\St_2V}\int_{\St_2V}g(\;R_{Y,\,X}X,\;Y\;)\;\vol_\FS(\,X,\,Y\,)
  $$
  Since the sectional curvature function $\sec_R$ is constant along the
  fibers of the projection $\pr_\Gr$ to the Grassmannian both integral
  formulas still hold true with $\St_2V$ replaced by $\Gr_2V$.
 \end{Theorem}

 \proof
 Recall that the projection $\St_2V\longrightarrow S_V,\,(X,Y)\longmapsto X,$
 is a Riemannian submersion for the Fubini--Study metric $g^\FS$ on the
 Stiefel manifold $\St_2V$ of orthonormal $2$--frames in $V$, hence we
 can replace the integral over $\St_2V$ by a double integral over
 spheres. Somewhat more precisely the inner integral is over the unit
 sphere $S_{\{\,X\,\}^\perp}$ of the $(m-1)$--dimensional orthogonal
 complement $\{\,X\,\}^\perp$ to the integration variable $X\,\in\,S_V$
 of the outer integral over $S_V$. Evidently inner and outer integration
 both integrate a homogeneous polynomial of degree $2k$ over a unit sphere
 so that Corollaries \ref{igauss} and \ref{sphint} together imply:
 \begin{eqnarray*}
  \lefteqn{\frac1{\Vol\;\St_2V}\int_{\St_2V}g(\;R_{Y,\,X}X,\;Y\;)^k
  \,\vol_\FS(\,X,\,Y\,)}\qquad
  &&
  \\
  &=&
  \frac1{\Vol\;S_V}\,\int_{S_V}\Big[\;\frac1{\Vol\;S_{\{\,X\,\}^\perp}}
  \,\int_{S_{\{\,X\,\}^\perp}}g(\;R_{Y,\,X}X,\;Y\;)^k\,\vol_g(\,Y\,)
  \;\Big]\,\vol_g(\,X\,)
  \\
  &=&
  \left.\frac{(\,-\,\Delta_g\,)^k}{4^k\,k!\,[\,\frac m2\,+\,k\,-\,1\,]_k}
  \right|_0\;\frac{k!}{2^k\,[\,\frac{m-1}2\,+\,k\,-\,1\,]_k}\;
  \\
  &&\qquad\left(\;X\,\longmapsto\,\mathrm{res}_{t\,=\,0}\Big[\;
  \frac{dt}{t^{k+1}}\,\exp\Big(\,\sum_{r>0}\frac{(2t)^r}{2r}\,
   \tr_V(\,R_{\,\cdot\,,\,X}X\,)^r\,\Big)\;\Big]\;\right)
 \end{eqnarray*}
 In this calculation the factor $2^k\,[\,\frac{m-1}2\,+\,k\,-\,1\,]_k$
 is the conversion factor of the crucial observation (\ref{conv}), which
 relates the integration over $S_{\{\,X\,\}^\perp}$ with the integration
 over $\{\,X\,\}^\perp$, the factor $k!$ on the other hand simply reflects
 that we are integrating $g(R_{Y,X}X,Y)^k$ over $Y$ instead of $\frac1{k!}
 \,g(R_{Y,X}X,Y)^k$ as in Corollary \ref{igauss}. In the resulting expression
 we have used a formal residue notation in order to pick up the coefficient
 of $t^k$ in the generating function:
 $$
  \mathrm{res}_{t\,=\,0}\Big[\;\frac{dt}{t^{k+1}}\,\exp\Big(\,
  \sum_{r>0}\frac{(2t)^r}{2r}\,\tr_V(\,R_{\,\cdot\,,\,X}X\,)^r\,\Big)\;\Big]
 $$
 Strictly speaking it is not necessary to single out the coefficient of
 $t^k$ in this way, because every $t^r$ comes along with a homogeneous
 polynomial in $X$ of degree $2r$ so that the coefficient of $t^k$ in the
 exponential is accompanied by a homogeneous polynomial in $X$ of degree
 $2k$. On the other hand the operator $\left.(-\Delta_g)^k \right|_0$ sees
 nothing else but homogeneous polynomials of degree $2k$ so that we may
 forget about picking up the coefficient of $t^k$ and may simply set
 $t\,=\,\frac12$ multiplying the overall result by $2^k$ in compensation.
 In this way we obtain eventually:
 \begin{eqnarray*}
  \lefteqn{\frac1{\Vol\;\St_2V}\int_{\St_2V}g(\;R_{Y,\,X}X,\;Y\;)^k
  \,\vol_\FS(\,X,\,Y\,)}\qquad
  &&
  \\
  &=&
  \left.\frac{(\,-\,\Delta_g\,)^k}{4^k\,[\,\frac m2\,+\,k\,-\,1\,]_k\,
  [\,\frac{m-1}2\,+\,k\,-\,1\,]_k}\right|_0\;\left(\;X\,\longmapsto\,
  \exp\Big(\,\sum_{r>0}\frac1{2r}\,\tr_V(\,R_{\,\cdot\,,\,X}X\,)^r\,\Big)
  \;\right)
  \\
  &=&
  \left.\frac{(\,-\,\Delta_g\,)^k}{[\,m\,+\,2k\,-\,2\,]_{2k}}\right|_0\;
  \left(\;X\,\longmapsto\,\exp\Big(\,\sum_{r>0}\frac1{2r}\,
  \tr_V(\,R_{\,\cdot\,,\,X}X\,)^r\,\Big)\;\right)
 \end{eqnarray*}
 In order to verify the statement about the scalar curvature of the algebraic
 curvature tensor $R$ we choose an orthonormal basis $X_1,\,\ldots,\,X_m$ of
 $V$ so that we may rewrite the trace of the Jacobi operator $\tr_V(\,R_{\,
 \cdot\,,\,X}X\,)\,=\,\sum g(\,R_{X_\mu,\,X}X,\,X_\mu\,)$ in the evaluation
 of the expression:
 $$
  (\,-\,\Delta_g\,)\left(\;X\,\longmapsto\,1\;+\;\frac12\,
  \sum_{\mu\,=\,1}^mg(\,R_{X_\mu,\,X}X,\,X_\mu\,)\;+\;\ldots\;\right)
  \;\;=\;\;
  \sum_{\mu,\,\nu\,=\,1}^mg(\,R_{X_\mu,\,X_\nu}X_\nu,\,X_\mu\,)
 $$
 \vskip-15pt
 \qed

 \pfill
 Taking up the lead given in the preceding proof we may define the
 generating power series of the moments of the sectional curvature
 of an algebraic curvature tensor $R$ on $V$ by:
 \begin{eqnarray*}
  \exp\Big(\,\sum_{r>0}\frac1{2r}\,\tr_V(\,R_{\,\cdot\,,\,X}X\,)^r\,\Big)
  &=&
  \exp\Big(\,-\,\frac12\,\sum_{r>0}\frac{(-1)^{r+1}}{r}
  \,\tr_V(\,-\,R_{\,\cdot\,,\,X}X\,)^r\,\Big)
  \\
  &=&
  \exp\Big(\;-\,\frac12\,\tr_V\;\log(\,\id\,-\,R_{\,\cdot\,,\,X}X\,)\;\Big)
 \end{eqnarray*}
 The standard matrix identity $\exp(\,\tr_VA\,)\,=\,\det_V(\,\exp\,A\,)$
 converts this definition into:

 \begin{Remark}[Generating Power Series for Sectional Curvature Moments]
 \hfill\label{gp}\break
  Consider an algebraic curvature tensor $R$ on a euclidean vector space
  $V$. In the sense of Theorem \ref{gins} the sectional curvature moments
  of $R$ have the following generating power series:
  $$
   \det^{-\,\frac12}_V\Big(\;\id\;-\;R_{\,\cdot\,,\,X}X\;\Big)
   \;\;:=\;\;
   \exp\Big(\,\sum_{r>0}\frac1{2r}\,\tr_V(\,R_{\,\cdot\,,\,X}X\,)^r\,\Big)
  $$
 \end{Remark}

 \pfill
 At the end of section we want to reformulate Theorem \ref{gins} in a way
 useful for the calculation of the moments of the sectional curvature of
 a compact Riemannian manifold $M$ of dimension $m$ with Riemannian metric
 $g$. The Grassmann bundle of $2$--planes in $TM$
 $$
  \Gr_2TM
  \;\;:=\;\;
  \{\;(\,p,\,\Sigma\,)\;|\;\;p\,\in\,M\textrm{\ and\ }\Sigma\,\in\,
  \Gr_2T_pM\;\;\}
 $$
 carries a natural Riemannian metric again denoted by $g$, which makes the
 canonical projection $\pr:\,\Gr_2TM\longrightarrow M$ a Riemannian submersion
 and induces the Fubini--Study metric $g^\FS$ on every fiber $\Gr_2T_pM$. The
 only thing we still have to specify is the orthogonal complement to the
 vertical tangent space, which we choose to be given by parallel transport
 of $2$--planes with respect to the Levi--Civita connection on $M$. In this
 way $\Gr_2TM$ becomes a compact Riemannian manifold endowed with a canonical
 function, namely the sectional curvature
 $$
  \sec:\quad\Gr_2TM\;\longrightarrow\;\R,\qquad
  (\,p,\,\Sigma\,)\;\longmapsto\;\sec_{R_p}(\;\Sigma\;)
 $$
 where $R_p$ is the curvature tensor of $M$ in the point $p$. In turn we may
 consider the sectional curvature of a compact Riemannian manifold $M$ as a
 random variable, whose moments and distribution form a legitimate object of
 study. The pinching constant for $M$ for example equals the quotient of the
 minimum and maximum value of $\sec$ provided $\sec_{R_p}(\Sigma)\,>\,0$ in
 at least one $2$--plane $(\,p,\,\Sigma\,)\,\in\,\Gr_2TM$.

 In this context Theorem \ref{gins} tells us that the moments of the
 sectional curvature of a Riemannian manifold $M$ can be calculated by
 integrating {\em natural} polynomials in the curvature tensor over $M$.
 Naturality in this context refers to the fact that these polynomials are
 defined for {\em all} euclidean vector spaces of dimension $m$ and are
 invariant under {\em all} isometries between these spaces. In particular
 natural polynomials are defined for all the tangent spaces $T_pM$ of
 Riemannian manifolds $M$ of dimension $m$, in turn their integrals over
 $M$ give rise to local Riemannian invariants for compact Riemannian
 manifolds in the sense of \cite{w}.

 \begin{Corollary}[Sectional Curvature Moments of Riemannian Manifolds]
 \hfill\label{mrm}\break
  For every euclidean vector space $V$ of dimension $m$ with scalar product
  $g$ there exist natural $\textbf{O}(\,V,\,g\,)$--invariant polynomials
  $\Psi^V_k$ of degree $k\,\in\,\N_0$ on the vector space $\Kr\,V$ of
  algebraic curvature tensors on $V$ such that the $k$--th moment of
  the sectional curvature equals
  $$
   \frac1{\Vol\;\Gr_2TM}\int_{\Gr_2TM}\sec_{R_p}^k(\;\Sigma\;)\;
   \vol_g(\,p,\,\Sigma\,)
   \;\;=\;\;
   \frac1{\Vol\;M}\int_M\Psi^{T_pM}_k(\;R_p\;)\,\vol_g(\,p\,)
  $$
  for every compact Riemannian manifold $M$ of dimension $m$, to wit
  the polynomial $\Psi^V_k$ reads:
  $$
   \Psi^V_k(\;R\;)
   \;\;:=\;\;
   \left.\frac{(\,-\,\Delta_g\,)^k}{[\,m\,+\,2k\,-\,2\,]_{2k}}\,\right|_0
   \left(\,X\,\longmapsto\,\exp\Big(\,\sum_{r>0}\frac1{2r}\,
   \tr_V(\,R_{\,\cdot\,,\,X}X\,)^r\,\Big)\,\right)
  $$
 \end{Corollary}
 
 \pfill
 Perhaps the most interesting aspect of this corollary is that it is
 possible, at least in principle, to calculate the pinching constant
 of a compact Riemannian manifold $M$ without calculating not even a
 single sectional curvature $\sec_{R_p}(\Sigma)$ explicitly. In fact
 it could turn out to be feasible in some examples to calculate all
 the higher moments of the random variable $\sec$ on $\Gr_2TM$, the
 minimum and maximum of the sectional curvature can then be reconstructed
 from these moments. The maximum of the absolute value of $\sec:\,\Gr_2TM
 \longrightarrow\R$ say is given by
 $$
  \max_{(\,p,\,\Sigma\,)\,\in\,\Gr_2TM}|\;\sec_{R_p}(\,\Sigma\,)\;|
  \;\;=\;\;
  \lim_{k\,\longrightarrow\,\infty}
  \sqrt[2k]{\frac1{\Vol\;\Gr_2TM}\int_{\Gr_2TM}\sec_{R_p}^{2k}(\;\Sigma\;)\;
   \vol_g(\,p,\,\Sigma\,)}
 $$
 and very similar formulas can be used to calculate the minimum and
 maximum of $\sec$.
\section{Examples of Sectional Curvature Distributions}
\label{ex}
 Explicit calculations in differential geometry are in general restricted
 to the relatively small class of Riemannian homogeneous spaces, in which
 all points look alike in the sense that the isometry group acts transitively.
 Homogeneous spaces are thus examples of choice to illustrate the calculation
 of the moments of sectional curvature discussed in Theorem \ref{gins}, not
 to the least because the eventual integration over the manifold in Corollary
 \ref{mrm} makes no difference at all. In this section we will calculate the
 moments and the distribution of the sectional curvature for all compact rank
 one symmetric spaces collectively known under the acronym CROSS: the round
 spheres $S^m$, the complex and quaternionic projective spaces $\C P^n$
 and $\H P^n$ with $n\,\geq\,2$ and the exceptional Cayley projective
 plane $\O P^2$. Moreover we will briefly discuss the problems arising
 in the case of Riemannian symmetric spaces of higher rank. A detailed
 introduction to the general theory of symmetric spaces including Cartan's
 classification can be found in the classical reference \cite{h}.
 
 \pfill
 Recall to begin with that the falling factorial polynomials
 $[\,x\,]_s\,:=\,x(x-1)\ldots(x-s+1)$ with parameter $s\,\in\,\Z$
 allow us to generalize the binomial coefficients via ${x\choose s}
 \,:=\,\frac1{s!}[x]_s$ to all $x\,\in\,\R$ and $s\,\in\,\Z$.
 One of the best--known of all the numerous identities satisfied
 by binomial coefficients extends directly to all $x,\,y\,\in\,\R$,
 simply because it depends polynomially on $x,\,y$:
 $$
  \sum_{s\,=\,0}^k{x\choose s}\,{y\choose k\,-\,s}
  \;\;=\;\;
  {x\,+\,y\choose k}
 $$
 In this form this classical binomial identity becomes important in the
 combinatorial argument
 \begin{eqnarray}
  \sum_{s\,=\,0}^k\frac{(-1)^s}{x\,+\,s}\,{k\choose s}
  &=&
  \frac{k!}{[\,x\,+\,k\,]_{k+1}}\,\sum_{s\,=\,0}^k(-1)^s\,
  \frac{[\,x\,+\,s\,-\,1\,]_s}{s!}\,\frac{[\,x\,+\,k\,]_{k-s}}{(k-s)!}
  \label{cx}
  \\
  &=&
  \frac{k!}{[\,x\,+\,k\,]_{k+1}}\,\sum_{s\,=\,0}^k
  {-x\choose s}\,{x\,+\,k\choose k-s}
  \;\;=\;\;
  \frac{k!}{[\,x\,+\,k\,]_{k+1}}
  \nonumber
 \end{eqnarray}
 valid for all $k\,\in\,\N_0$, which allows us to evaluate the following
 integral for all $a,\,b,\,k\,\in\,\N_0$:
 \begin{eqnarray*}
  \int_0^1(\,u^2\,)^a\,(\,1-u^2\,)^b\,(\,1+3u^2\,)^k\,du
  &=&
  \sum_{r\,=\,0}^k4^r\,{k\choose r}
  \int_0^1(\,u^2\,)^{a+r}\,(\,1-u^2\,)^{b+k-r}\,du
  \\
  &=&
  \sum_{r\,=\,0}^k4^r\,{k\choose r}\,\left(\sum_{s\,=\,0}^{b+k-r}
  \frac{(-1)^s}{2a+2r+2s+1}\,{b+k-r\choose s}\right)
  \\
  &=&
  \frac12\sum_{r\,=\,0}^k4^r\,{k\choose r}\,
  \frac{[\,b\,+\,k\,-\,r\,]_{k-r}\,[\,b\,]_b\,[\,a\,+\,r\,-\,\frac12\,]_r}
  {[\,a\,+\,b\,+\,k\,+\,\frac12\,]_k\,[\,a\,+\,b\,+\,\frac12\,]_{b+1}}
  \\
  &=&
  \frac1{(\,2a\,+\,1\,)\,{a\,+\,b\,+\,\frac12\choose b}}
  \,\sum_{r\,=\,0}^k4^k\,\frac{{-\,\frac{2a\,+\,1}2\choose r}
  \,{-\,\frac{2b\,+\,2}2\choose k-r}}{{-\,\frac{2a\,+\,2b\,+\,3}2\choose k}}
 \end{eqnarray*}
 In the last line we have used $[\,x\,+\,s\,]_s\,=\,(-1)^s\,
 [\,-\frac{2x\,+\,2}2\,]_s$ thrice in order to simplify the result.
 Making the substitution $s\,=\,1+3u^2,\,ds\,=\,6\,u\,du,$ in the
 integral we conclude:
 \begin{equation}\label{prob}
  \frac{(\,2a\,+\,1\,)}6\,{a\,+\,b\,+\,\frac12\choose b}
  \,\int_1^4\sqrt{\frac{s-1}3}^{2a-1}\,\Big(\frac{4-s}3\Big)^b\,s^k\,ds
  \;\;=\;\;
  \sum_{r\,=\,0}^k4^k\,\frac{{-\,\frac{2a\,+\,1}2\choose r}
  \,{-\,\frac{2b\,+\,2}2\choose k-r}}{{-\,\frac{2a\,+\,2b\,+\,3}2\choose k}}
 \end{equation}
 En nuce this formula can be seen as a formula for the $k$--th order moments
 of a family of probability measures on the interval $[\,1,\,4\,]$ parametrized
 by $a,\,b\,\in\,\N_0$.

 \pfill
 The distribution of the sectional curvature on the round spheres $S^m$ of
 dimension $m$ and scalar curvature $\kappa\,>\,0$, the simplest examples
 of compact symmetric spaces of rank $1$, is too simple to be of any interest,
 because $S^m$ has constant sectional curvature $\sec\,\equiv\,\frac{\kappa}
 {m\,(m-1)}$. Nevertheless it may be instructive to apply Lemma \ref{gins}
 in the special case $\kappa\,=\,m(m-1)$ corresponding the round metric on
 $S^m$ of constant sectional curvature $\sec\,\equiv\,1$ to verify that
 all higher moments turn out to be equal to $1$. Instead of using the
 algebraic curvature tensor $R^{S^m}$ of the sphere $S^m$ of radius $1$
 directly we consider the so--called Jacobi operator
 $$
  R^{S^m}_{\,\cdot\,,\,X}X:\quad V\;\longrightarrow\;V,\qquad
  Y\;\longmapsto\;R^{S^m}_{Y,\,X}X\;\;=\;\;g(\,X,\,X\,)\,Y\;-\;g(\,X,\,Y\,)\,X
 $$
 for a fixed vector $X\,\in\,V$, which evidently has eigenvalues $0$ and
 $g(X,X)$ on the eigenspaces $\R\,X$ and $\{\,X\,\}^\perp$ respectively.
 In turn the generating power series of Remark \ref{gp} reads
 \begin{equation}\label{sm}
  \det_V^{-\,\frac12}\Big(\;\id\;-\;R^{S^m}_{\,\cdot\,,\,X}X\;\Big)
  \;\;=\;\;
  (\;1\;-\;g(\,X,\,X\,)\;)^{-\frac{m-1}2}
  \;\;=\;\;
  \sum_{k\,\geq\,0}{-\frac{m-1}2\choose k}(\,-\,g(\,X,\,X\,)\,)^k
 \end{equation}
 using once again Newton's power series expansion $(1+t)^x\,=\,
 \sum_{k\,\geq\,0}{x\choose k}\,t^k$. On the other hand the polynomial
 $X\longmapsto g(X,X)^k$ equals $1$ on the unit sphere $S_V\,\subset\,V$
 of a euclidean vector space $V$ of dimension $m$ and thus integrates over
 $S_V$ to $\Vol\,S_V$ for all $k\,\in\,\N_0$. Reading Corollary \ref{sphint}
 in light of this trivial observation we obtain directly:
 \begin{equation}\label{scalar}
  \left.\frac{(\,-\,\Delta_g\,)^k}{[\,m\,+\,2k\,-\,2\,]_{2k}}\right|_0
  \Big(\;X\,\longmapsto\,g(\,X,\,X\,)^k\;\Big)
  \;\;=\;\;
  \frac{4^k\,k!\,[\,\frac m2\,+\,k\,-\,1\,]_k}
  {[\,m\,+\,2k\,-\,2\,]_{2k}}
  \;\;=\;\;
  \frac{(\,-\,1\,)^k}{{-\,\frac{m-1}2\choose k}}
 \end{equation}
 Comparing this result with equation (\ref{sm}) for the generating power
 series of sectional curvature moments we conclude that all moments of
 the round metric on $S^m$ equal $1$.  

 Turning from round spheres to more complicated examples we recall that the
 complex projective spaces $\C P^n,\,n\,\geq\,2,$ of dimensions $m\,=\,2n$
 are the prototypical examples of K\"ahler manifolds. In particular all their
 tangent spaces $T_p\C P^n$ carry an orthogonal complex structure $I\,\in\,
 \End\,T_p\C P^n$ satisfying $I^2\,=\,-\id$ as well as $g(IX,IY)\,=\,g(X,Y)$.
 According to \cite{bes} the algebraic curvature tensor of the unique symmetric
 Riemannian metric $g$ on $\C P^n$ of scalar curvature $\kappa\,>\,0$ can
 be written in terms of this orthogonal complex structure:
 \begin{eqnarray*}
  \lefteqn{R^{\C P^n}_{X,\,Y}Z}
  &&
  \\
  &\!=\!&
  -\,\frac\kappa{4\,n\,(n+1)}\Big(\,g(X,Z)Y\,-\,g(Y,Z)X\,+\,g(IX,Z)IY
  \,-\,g(IY,Z)IX\,+\,2g(IX,Y)IZ\,\Big)
 \end{eqnarray*}
 In consequence the Jacobi operators $R_{\,\cdot\,,\,X}X:\,T_p\C P^n
 \longrightarrow T_p\C P^n$ read for all $X\,\in\,T_p\C P^n$:
 $$
  Y\;\longmapsto\;R^{\C P^n}_{Y,\,X}X
  \;\;=\;\;
  -\,\frac\kappa{4\,n\,(n+1)}\,\Big(\,g(Y,X)X\,-\,g(X,X)Y
  \,-\,3\,g(Y,IX)IX\,\Big)
 $$
 With $\C P^n$ being a rank $1$ symmetric space the spectrum of the
 Jacobi operator $R^{\C P^n}_{\,\cdot\,,\,X}X$ depends of course only
 on the norm square $g(X,X)$ of the argument vector $X\,\in\,T_p\C P^n$.
 Specifically we find for the scalar curvature $\kappa\,=\,4\,n\,(n+1)$ the
 eigenvalues $0$, $4g(X,X)$ and $g(X,X)$ on the eigenspaces $\R\,X$, $\R\,IX$
 and $\{\,X,\,IX\,\}^\perp$ respectively so that the sectional curvature
 for this value of $\kappa$ takes values in the interval $[\,1,\,4\,]$.
 The generating power series of the moments of the sectional curvature
 introduced in Remark \ref{gp} thus reads:
 \begin{eqnarray*}
  \exp\Big(\,\sum_{r>0}\frac1{2r}\,
  \tr_V(\,R^{\C P^n}_{\,\cdot\,,\,X}X\,)^r\,\Big)
  &=&
  (\,1\,-\,4\,g(\,X,\,X\,)\,)^{-\frac12}
  \,(\,1\,-\,g(\,X,\,X\,)\,)^{-\frac{2n-2}2}
  \\
  &=&
  \sum_{k\,\geq\,0}\left[\;\sum_{r\,=\,0}^k4^r\,{-\frac12\choose r}
  \,{-\frac{2n-2}2\choose k-r}\;\right]\,(\;-\,g(X,X)\;)^k
 \end{eqnarray*}
 Equation (\ref{scalar}) and Theorem \ref{gins} tell us the $k$--th order
 moment of the sectional curvature
 \begin{equation}\label{cpm}
  \frac1{\Vol\;\Gr_2T\C P^n}
  \int_{\Gr_2T\C P^n}\!\!\!\!\sec^k\;\vol_g
  \;\;=\;\;
  \sum_{r\,=\,0}^k4^r\,\frac{{-\,\frac12\choose r}
  \,{-\,\frac{2n-2}2\choose k-r}}{{-\,\frac{2n-1}2\choose k}}
 \end{equation}
 of $\C P^n$. Comparing this result with the explicit moments of the
 probability measures defined in equation (\ref{prob}) we find a match
 for $a\,=\,0$, $b\,=\,n-2$, the corresponding probability measure thus
 describes the distribution of the sectional curvature on $\C P^n$.

 Virtually nothing in this argument has to be changed to calculate the
 distribution of the sectional curvature on the quaternionic projective
 spaces $\H P^n$ of dimension $m\,=\,4n$. Instead of one orthogonal complex
 structure $I$ the tangent spaces $T_p\H P^n$ of $\H P^n$ feature three such
 structures $I,\,J,\,K\,\in\,\End\,T_p\H P^n$ satisfying $IJ\,=\,K$, the
 only difference to $\C P^n$ is that $I,\,J,\,K$ are neither globally
 defined nor do we have a means to distinguished them. Nevertheless
 the algebraic curvature tensor of the symmetric metric $g$ on $\H P^n$
 of scalar curvature $\kappa\,>\,0$ is simply the extension of $R^{\C P^n}$
 to three orthogonal complex structures
 \begin{eqnarray*}
  R^{\H P^n}_{X,\,Y}Z
  &=&
  -\frac\kappa{16\,n\,(n+2)}\,\Big(\;g(\,X,\,Z)\,Y\;-\;g(\,Y,\,Z)\,X
  \\[-5pt]
  &&
  \hskip68pt
  \;+\;g(\,\hbox to11pt{\hfil$I$\hfil}X,\,Z\,)\,\hbox to11pt{\hfil$I$\hfil}Y
  -\;g(\,\hbox to11pt{\hfil$I$\hfil}Y,\,Z\,)\,\hbox to11pt{\hfil$I$\hfil}X
  +\;2\,g(\,\hbox to11pt{\hfil$I$\hfil}X,\,Y\,)\,\hbox to11pt{\hfil$I$\hfil}Z
  \\
  &&
  \hskip68pt
  \;+\;g(\,\hbox to11pt{\hfil$J$\hfil}X,\,Z\,)\,\hbox to11pt{\hfil$J$\hfil}Y
  -\;g(\,\hbox to11pt{\hfil$J$\hfil}Y,\,Z\,)\,\hbox to11pt{\hfil$J$\hfil}X
  +\;2\,g(\,\hbox to11pt{\hfil$J$\hfil}X,\,Y\,)\,\hbox to11pt{\hfil$J$\hfil}Z
  \\[-3pt]
  &&
  \hskip68pt
  \;+\;g(\,\hbox to11pt{\hfil$K$}X,\,Z\,)\,\hbox to11pt{\hfil$K$}Y
  -\;g(\,\hbox to11pt{\hfil$K$}Y,\,Z\,)\,\hbox to11pt{\hfil$K$}X
  +\;2\,g(\,\hbox to11pt{\hfil$K$}X,\,Y\,)\,\hbox to11pt{\hfil$K$}Z\;\Big)
 \end{eqnarray*}
 compare \cite{bes}. Similar to the case of $\C P^n$ the Jacobi operator
 $R^{\H P^n}_{\,\cdot\,,\,X}X$ has eigenvalues $0$, $4g(X,X)$ and $g(X,X)$
 on $\R\,X$, $\span\,\{\,IX,\,JX,\,KX\,\}$ and $\{\,X,\,IX,\,JX,\,KX\,\}^\perp$
 respectively for the special value $\kappa\,=\,16\,n\,(n+2)$ of the scalar
 curvature, for which the sectional curvature takes values in the interval
 $[\,1,\,4\,]$. In consequence we find the generating function
 $$
  \exp\Big(\,\sum_{r>0}\frac1{2r}\,
  \tr_V(\,R^{\H P^n}_{\,\cdot\,,\,X}X\,)^r\,\Big)
  \;\;=\;\;
  (\,1\,-\,4\,g(\,X,\,X\,)\,)^{-\frac32}
  \,(\,1\,-\,g(\,X,\,X\,)\,)^{-\frac{4n-4}2}
 $$
 and conclude for the $k$--th order moment of the sectional
 curvature of $\H P^n$:
 \begin{equation}\label{hpm}
  \frac1{\Vol\;\Gr_2T\H P^n}
  \int_{\Gr_2T\H P^n}\!\!\!\!\sec^k\;\vol_g
  \;\;=\;\;
  \sum_{r\,=\,0}^k4^r\,\frac{{-\,\frac32\choose r}
  \,{-\,\frac{4n-4}2\choose k-r}}{{-\,\frac{4n-1}2\choose k}}
 \end{equation}
 In terms of the probability measures defined in equation (\ref{prob})
 these moments correspond to the probability measure with parameters
 $a\,=\,1$ and $b\,=\,2n-3$.

 The curvature tensor of the only remaining compact symmetric space of rank
 $1$, the Cayley projective plane $\O P^2$ of dimension $m\,=\,16$, is more
 difficult to describe, because we would need to discuss some properties of
 spinors in dimension $9$ first. For the unique symmetric Riemannian metric
 $g$ on $\O P^2$ of scalar curvature $\kappa\,=\,576$ the Jacobi operators
 $R^{\O P^2}_{\,\cdot\,,\,X}X$ have nevertheless the eigenvalues $0$,
 $4g(X,X)$ and $g(X,X)$ of multiplicities $1$, $7$ and $8$ respectively.
 The generating function of sectional curvature moments with these
 multiplicities
 $$
  \exp\Big(\,\sum_{r>0}\frac1{2r}\,
  \tr_V(\,R^{\O P^2}_{\,\cdot\,,\,X}X\,)^r\,\Big)
  \;\;=\;\;
  (\,1\,-\,4\,g(\,X,\,X\,)\,)^{-\frac72}\,(\,1\,-\,g(\,X,\,X\,)\,)^{-\frac82}
 $$
 corresponds to the probability measure with parameters $a\,=\,b\,=\,3$
 and the moments:
 \begin{equation}\label{opm}
  \frac1{\Vol\;\Gr_2T\O P^2}
  \int_{\Gr_2T\O P^2}\!\!\!\!\sec^k\;\vol_g
  \;\;=\;\;
  \sum_{r\,=\,0}^k4^r\,\frac{{-\,\frac72\choose r}
  \,{-\,\frac82\choose k-r}}{{-\,\frac{15}2\choose k}}
 \end{equation}

 \begin{Corollary}[Examples of Sectional Curvature Densities]
 \hfill\label{cross}\break
  The sectional curvature $\sec:\,\Gr_2TM\longrightarrow[\,\min,\,\max\,]$
  of a Riemannian metric $g$ on a compact Riemannian manifold $M$ can be
  thought of as a random variable with associated probability measure on
  $[\,\min,\,\max\,]$. In this way the unique symmetric Riemannian metrics
  with sectional curvatures in $[\,1,\,4\,]$ on the complex and quaternionic
  projective spaces $\C P^n$ and $\H P^n$ with $n\geq2$ as well as the Cayley
  projective plane $\mathbb{O}P^2$ define the probability measures:
  \begin{eqnarray*}
   \frac1{\Vol\;\Gr_2T\C P^n}
   \int_{\Gr_2T\C P^n}\!\!\!\!F(\,\sec\,)\;\vol_g
   &=&
   \hbox to 64pt{${\displaystyle\frac16\;\,{n-\frac32\choose n-2}}$\hfill}
   \int_1^4\hbox to50pt{$\displaystyle\sqrt{\frac{s-1}3}^{-1}$\hss}
   \hbox to64pt{$\displaystyle\Big(\frac{4-s}3\Big)^{n-2}$\hfil}F(\,s\,)\;ds
   \\[2pt]
   \frac1{\Vol\;\Gr_2T\H P^n}
   \int_{\Gr_2T\H P^n}\!\!\!\!F(\,\sec\,)\;\vol_g
   &=&
   \hbox to 64pt{${\displaystyle\frac36\;{2n-\frac32\choose 2n-3}}$\hfill}
   \int_1^4\hbox to50pt{$\displaystyle\sqrt{\frac{s-1}3}$\hfil}
   \hbox to64pt{$\displaystyle\Big(\frac{4-s}3\Big)^{2n-3}$\hfil}F(\,s\,)\;ds
   \\[2pt]
   \frac1{\Vol\;\Gr_2T\O P^2}
   \int_{\Gr_2T\O P^2}\!\!\!\!F(\,\sec\,)\;\vol_g
   &=&
   \hbox to 64pt{${\displaystyle\frac76\quad\;{\frac{13}2\choose 3}}$\hfill}
   \int_1^4\hbox to50pt{$\displaystyle\sqrt{\frac{s-1}3}^5$\hfil}
   \hbox to64pt{$\displaystyle\Big(\frac{4-s}3\Big)^3$\hfil}F(\,s\,)\;ds
  \end{eqnarray*}
 \end{Corollary}

 \pfill
 The calculation of the sectional curvature densities of compact Riemannian
 symmetric spaces of higher rank is significantly more involved, however we
 want to discuss the general setup at the end of this section. In general every
 Killing field $X$ for a Riemannian metric $g$ on a connected manifold $M$ is
 automatically an affine Killing field for the Levi--Civita connection
 $$
  0
  \;\;=\;\;
  (\,\mathfrak{Lie}_X\nabla\,)_YZ
  \;\;:=\;\;
  [\,X,\,\nabla_YZ\,]\;-\;\nabla_{[\,X,\,Y\,]}Z\;-\;\nabla_Y[\,X,\,Z\,]
  \;\;=\;\;
  \nabla^2_{Y,\,Z}X\;+\;R_{X,Y}Z
 $$
 and thus satisfies the extended Killing equation with the auxiliary
 section $\X\,:=\,\nabla X$
 $$
  \nabla_YX\;\;=\;\;\X\,Y
  \qquad\qquad
  \nabla_Y\X\;\;=\;\;-\,R_{X,\,Y}
 $$
 of the bundle $\End_{\mathrm{skew}}TM\,\subset\,\End\,TM$ of skew--symmetric
 endomorphisms of $TM$. With $\nabla$ being torsion free the Lie bracket of
 two extended Killing fields $\X\oplus X$ and $\Y\oplus Y$ reads
 \begin{equation}\label{br}
  [\;\X\,\oplus\,X,\;\Y\,\oplus\,Y\;]
  \;\;=\;\;
  \Big(\;R_{X,Y}\;-\;[\,\X,\,\Y\,]\;\Big)\,\oplus\,
  \Big(\;\Y\,X\;-\;\X\,Y\;\Big)
 \end{equation}
 where $[\,\X,\,\Y\,]$ is simply the pointwise commutator of endomorphisms,
 because
 $$
  \nabla_Z[\,X,\,Y\,]
  \;\;=\;\;
  \nabla_Z(\;\Y\,X\,-\,\X\,Y\;)
  \;\;=\;\;
  -\;R_{Y,\,Z}X\;+\;\Y\,(\,\X\,Z\,)\;+\;R_{X,\,Z}Y\;-\;\X\,(\,\Y\,Z\,)
 $$
 for every vector field $Z\,\in\,\Gamma(\,TM\,)$. In other words the vector
 bundle $\End_{\mathrm{skew}}TM\,\oplus\,TM$ becomes a bundle of algebras,
 but not in general a bundle of Lie algebras, under the pointwise bracket
 defined in equation (\ref{br}) with the property that the evaluation in
 $p\,\in\,M$
 $$
  \iota_p:\quad\g\;\longmapsto\;\End_{\mathrm{skew}}T_pM\,\oplus\,T_pM,
  \qquad X\;\longmapsto\;(\,\nabla\,X\,)_p\,\oplus\,X_p
 $$
 is an injective {\em antihomomorphism} of algebras, compare \cite{bgv}.
 For Riemannian homogeneous spaces the composition of $\iota_p$ with the
 projection to $T_pM$ is surjective, what makes symmetric spaces rather
 special among homogeneous spaces however is that the image of $\g$ is
 compatible
 $$
  \g
  \;\;\cong\;\;
  \iota_p(\;\g\;)
  \;\;=\;\;
  \Big(\;\iota_p(\,\g\,)\,\cap\,\End_{\mathrm{skew}}T_pM\;\Big)
  \,\oplus\,\Big(\;\iota_p(\,\g\,)\,\cap\,T_pM\;\Big)
 $$
 with the direct sum for all $p\,\in\,M$ and thus invariant under the
 involutive automorphism $\theta_p$ of the algebra $\End_{\mathrm{skew}}T_pM
 \oplus T_pM$ defined by $\theta_p(\,\X\oplus X\,)\,:=\,\X\oplus(-X)$. In
 consequence the Lie algebra $\g$ of Killing vector fields of a connected
 Riemannian symmetric space $M$ decomposes as a vector space into the direct
 sum $\g\,=\,\k_p\oplus\p_p$ of the preimages $\k_p$ and $\p_p$ of
 $\End_{\mathrm{skew}}T_pM$ and $T_pM$ respectively under $\iota_p$,
 which agree with the $(+1)$ and $(-1)$--eigenspaces respectively of
 the involutive automorphism $\theta_p$ of $\End_{\mathrm{skew}}T_pM
 \oplus T_pM$ restricted to $\g$.

 For the purpose of calculating moments of sectional curvature the most
 interesting consequence of the preceding discussion of symmetric spaces
 is certainly that the sectional curvature of a Riemannian symmetric space
 is completely determined by the algebra structure of its Lie algebra $\g$
 of Killing vector fields. More precisely there exists a unique invariant
 scalar product on $\g$ such that the evaluation map $\mathrm{ev}_p:\,\p_p
 \stackrel\cong\longrightarrow T_pM,\,X\longmapsto X_p,$ is an isometry for
 some and hence all $p\,\in\,M$, moreover this isometry identifies the
 curvature tensor $R^{G/K}_p$ in a point $p\,\in\,M$ with the so--called
 triple product \cite{h} defined on the subspace $\p_p\,\subset\,\g$
 $$
  (\,\ad^2X\,)\,Y
  \;\;=\;\;
  [\;0\,\oplus\,X,\;[\;0\,\oplus\,X,\;0\,\oplus\,Y\;]\;]
  \;\;=\;\;
  0\,\oplus\,(\,R^{G/K}_p\,)_{X,\,Y}X
 $$
 for all $X,\,Y\,\in\,\p_p$. In particular the generating power series
 of the sectional curvature moments defined in Remark \ref{gp} can be
 replaced by $\det^{-\frac12}_{\p_p}(\,\id\,+\,\ad^2X\,)$ in the algebraic
 model $\p_p\,\cong\,T_pM$ so that Theorem \ref{gins} implies the following
 generating formal power series
 $$
  \sum_{k\,\geq\,0}(\,-t\,)^k \;{-\,\frac{m-1}2\choose k}
  \;\Psi_k(\;R^{G/K}_p\;)
  \;\;=\;\;
  \frac1{\Vol\,S_\p}\,\int_{S_\p}
  \det^{-\frac12}_\p(\;\id\,+\,t\,\ad^2X\;)\,\vol_g(\,X\,)
 $$
 for the moments of the sectional curvature of the algebraic curvature
 tensor $R^{G/K}_p$ corresponding to the Riemannian symmetric space $M\,=\,
 G/K$. Combining this explicit description of the $k$--th order moments
 $\Psi_k(\,R^{G/K}_p\,)$ of the sectional curvature of $M$ with a suitable
 variant of the classical Integration Formula of Weyl for compact Lie groups
 \cite{bgv} we obtain eventually the following formula for the moments of
 the sectional curvature of $M\,=\,G/K$:

 \begin{Lemma}[Weyl's Integration Formula for Sectional Curvature Moments]
 \hfill\label{wif}\break
  Consider a compact Riemannian symmetric space $M\,=\,G/K$ of dimension $m$
  with isometry group $G\,:=\,\mathrm{Isom}\,M$ and stabilizer $K\,\subset\,G$
  of a given base point $p\,\in\,M$. The Lie algebra of Killing vector fields
  decomposes orthogonally $\g\,=\,\k\oplus\p$ with respect to the unique
  invariant scalar product $g^{\mathrm{ext}}$ on $\g$, which makes the
  evaluation $\mathrm{ev}_p:\,\p\longrightarrow T_pM$ in the base point
  an isometry. Choosing a maximal abelian subalgebra $\a\,\subset\,\p$
  with centralizer and normalizer $\mathrm{Zent}\,\a\,\subset\,\mathrm{Norm}
  \,\a\,\subset\,K$ we may encode the sectional curvature moments of $M$ in
  the series
  \begin{eqnarray*}
   \lefteqn{\sum_{k\,\geq\,0}(\,-t\,)^k\;{-\,\frac{m-1}2\choose k}\;
    \Psi_k(\;R^{G/K}_p\;)}\quad
   &&
   \\
   &=&
   \frac1{|\,\mathfrak{W}_{\a}\,|}
   \frac{\Vol\,(\,K,\,g^{\mathrm{ext}\,)}}
   {\Vol\,S_\p\,\cdot\,\Vol\,(\,\mathrm{Zent}\,\a,\,g^{\mathrm{ext}}\,)}
   \int_{S_\a}\det^{-\frac12}_\p(\;\id\,+\,t\,\ad^2X\;)\,
   \sqrt{\det^*_\p(\;-\,\ad^2X\;)}\,\vol_g(X)
  \end{eqnarray*}
  where $\mathfrak{W}_\a\,:=\,\mathrm{Norm}\,\a/\mathrm{Zent}\,\a$ is the Weyl
  group associated to $\a\,\subset\,\p$ and $\det^*_\p(-\ad^2X)$ is the product
  of all non--zero and hence positive eigenvalues of $-\ad^2X$ on $\p$.
 \end{Lemma}

 \pfill
 In contrast to the standard Integration Formula for compact Lie groups the
 square root of the orbit volume ratio $\det^*_\p(-\ad^2X)$ in the modified
 Integration Formula for compact Riemannian symmetric spaces can not be
 replaced by a polynomial in $X$, because the integrand does not change sign
 under the Weyl group. More than a mere annoyance this unexpected twist is
 a serious obstacle in the explicit calculation of the moments of the sectional
 curvature, because the integration trick of Lemma \ref{sphint} does no longer
 apply. At least for symmetric spaces of rank $\dim\,\a\,=\,2$ this problem
 can be circumvented by replacing the integration over $S_\a$ by an integration
 over the circular arc $S_\a/\mathfrak{W}_\a$. In this way one obtains for
 example the moments
 $$
  \Psi_k(\,R^{\Gr_2\R^n}_p\,)
  \;\;=\;\;
  \frac1{{-\frac{2n-5}2\choose k}}\,\sum_{{{\scriptstyle\mu,\,\nu\,\geq\,0}
  \atop{\scriptstyle\mu+\nu\,\leq\,\frac k2}}}\!\!\!
  \,\frac{(-1)^\nu}{4^\mu}\,\frac{n\,-\,3}{n\,+\,2\mu\,+\,2\nu\,-\,3}
  \,{-\frac{n-4}2\choose\mu}\,{-\frac12\choose\nu}
  \,{-\frac{n+2\mu+4\nu-2}2\choose k-2\mu-2\nu}
 $$
 of the sectional curvature for the Fubini--Study metric on the Grassmannians
 of $2$--planes in $\R^n$ with sectional curvatures in the interval
 $[\,0,\,2\,]$ and scalar curvature $\kappa\,=\,2(n-2)^2$, somewhat
 more complicated formulas hold true for the Grassmannians of $2$--planes
 in $\C^n$ and $\H^n$. For the Grassmannian $\Gr_2\R^4$ the explicit formula
 can be verified directly, because $S^2\times S^2$ endowed with half the
 product metric is a $2$--fold Riemannian cover of $\Gr_2\R^4$.
\section{Moments and Densities of Riemannian Products}
\label{mdp}
 Certainly the simplest way to produce new examples of Riemannian
 manifolds is to take products, hence it seems prudent to study this
 construction with a view on the moments of the sectional curvature
 in order to extend our stock of examples. Somewhat surprisingly it
 turns out that both the moments and the distribution of the sectional
 curvature behave very nicely under Riemannian products to the extent
 that there exists a simple probabilistic model for Riemannian products.
 The explicitly known values of the moments of sectional curvature for
 the examples $\C P^2,\,S^4$ together with the values for the products
 $S^2\times S^2$ and $S^1\times S^3$ calculated in this section will be
 used in due course to prove a weak version of a well--known theorem of
 Hitchin--Thorpe \cite{th} concerning $4$--dimensional Einstein manifolds.

 \pfill
 Recall that the canonical projections $\pr_M:\,M\times N\longrightarrow M,
 \,(\,p,\,q\,)\longmapsto p,$ and $\pr_N$ from the Cartesian product
 $M\times N$ of two Riemannian manifolds $M$ and $N$ with Riemannian
 metrics $g$ and $h$ respectively to its factors induces a canonical
 isomorphism of tangent spaces 
 $$
  T_{(\,p,\,q\,)}(\;M\,\times\,N\;)
  \;\stackrel\cong\longrightarrow\;T_pM\,\oplus\,T_qN,
  \qquad\left.\frac d{dt}\right|_0(\,p_t,\,q_t\,)\;\longmapsto\;
  \left.\frac d{dt}\right|_0p_t\;\oplus\;\left.\frac d{dt}\right|_0q_t
 $$
 which turns the Cartesian into the Riemannian product $M\times N$
 endowed with the metric:
 $$
  (\,g^M\,\oplus\,g^N\,)_{(\,p,\,q\,)}(\;X,\;Y\;)
  \;\;:=\;\;
  g^M_p(\;\pr_MX,\;\pr_MY\;)\;+\;g^N_q(\;\pr_NX,\;\pr_NY\;)
 $$
 The curvature $R^{M\times N}_{(\,p,\,q\,)}$ of the Levi--Civita connection
 $\nabla$ associated to this product metric in a point $(p,q)\,\in\,M\times N$
 decomposes likewise into the direct sum $R^M_p\,\oplus R^N_q$ in the sense
 $$
  (\,R^{M\times N}_{(\,p,\,q\,)}\,)_{X,\,Y}Z
  \;\;=\;\;
  (\,R^M_p\,)_{\pr_MX,\,\pr_MY}\pr_MZ\,\oplus\,
  (\,R^N_p\,)_{\pr_NX,\,\pr_NY}\pr_NZ
 $$
 for all $X,\,Y,\,Z\,\in\,T_{(p,q)}(M\times N)$. In turn the Jacobi operator
 associated to a tangent vector $X\,\in\,T_{(p,q)}(M\times N)$ reduces to the
 direct sum $R^{M\times N}_{\,\cdot\,,\,X}X:\,T_pM\oplus T_qN\longrightarrow
 T_pM\oplus T_qN$ of the Jacobi operators of the projections $\pr_MX\,\in\,
 T_pM$ and $\pr_NX\,\in\,T_qN$ so that the generating power series of
 sectional curvature moments of $M\times N$ introduced in Remark \ref{gp}
 \begin{eqnarray*}
  \lefteqn{\det^{-\frac12}_{T_pM\,\oplus\,T_qN}
  \Big(\;\id\;-\;R^{M\,\times\,N}_{\;\cdot\;,\,X}X\;\Big)}\qquad&&
  \\
  &=&
  \det^{-\frac12}_{T_pM}\Big(\;\id\;-\;R^M_{\;\cdot\;,\,\pr_M\,X}
  \pr_M\,X\;\Big)\;\cdot\;\det^{-\frac12}_{T_qN}
  \Big(\;\id\;-\;R^N_{\;\cdot\;,\,\pr_N\,X}\pr_N\,X\;\Big)
 \end{eqnarray*}
 becomes simply the product. As a direct consequence we obtain from Theorem
 \ref{gins}:
  
 \begin{Lemma}[Sectional Curvature Moments of Products]
 \hfill\label{mp}\break
  Let $M$ and $N$ be Riemannian manifolds of dimensions $m,\,n\,\geq\,2$
  respectively. The curvature tensor of the Riemannian product manifold
  $M\times N$ in a point $(p,q)\,\in\,M\times N$ decomposes into the direct
  sum $R^{M\times N}_{(\,p,\,q\,)}\,=\,R^M_p\,\oplus\,R^N_q$ of the curvature
  tensors of $M$ in $p$ and $N$ in $q$. Hence the pointwise moments of the
  sectional curvatures of $R^{M\times N}$ can be calculated via
  $$
   \Psi_k(\,R^{M\times N}_{(\,p,\,q\,)}\,)
   \;\;=\;\;
   \sum_{r\,=\,0}^k{k\choose r}\,\frac{[\,m+2r-2\,]_{2r}
   \,[\,n+2(k-r)-2\,]_{2(k-r)}}{[\,m+n+2k-2\,]_{2k}}
   \,\Psi_r(\,R^M_p\,)\,\Psi_{k-r}(\,R^N_q\,)
  $$
  and this formula still holds true for the integrated pointwise moments of
  Corollary \ref{mrm}.
 \end{Lemma}

 \pfill
 From the point of view of combinatorics the formula for the sectional
 curvature moments of Riemannian products is rather strange, however it
 turns out that this formula has a very neat model in probability theory.
 In order to discuss this model we recall a simple integration trick
 for integrating polynomials over the standard simplex in dimension
 $n\,\in\,\N_0$:
 $$
  \Delta_n
  \;\;:=\;\;
  \{\;\,(\,x_1,\,\ldots,\,x_n\,)\,\in\,\R^n\;\;|
  \;\;x_1,\,\ldots,\,x_n\,\geq\,0\textrm{\ and\ }x_1+\ldots+x_n\,\leq\,1\;\,\}
 $$
 Every polynomial in the coordinates $x_1,\,\ldots,\,x_n$ and $x_0\,:=\,
 1-x_1-\ldots-x_n$ integrates to
 \begin{equation}\label{if}
  \int_{\Delta_n}x_0^{k_0}\,x_1^{k_1}\,x_2^{k_2}\,\ldots\,x_n^{k_n}
  \;\;dx_1\,dx_2\,\ldots\,dx_n
  \;\;=\;\;
  \frac{k_0!\,k_1!\;k_2!\;\ldots\;k_n!}{(\,n+k_0+k_1+k_2+\ldots+k_n\,)!}
 \end{equation}
 against the Lebesgue measure on the standard simplex $\Delta_n\,
 \subset\,\R^n$. Needless to say this integration trick can be proved
 by straightforward induction on $n$ with the combinatorial argument
 (\ref{cx}) coming in handily in both the base and the induction step.
 In consequence the following density defines a probability measure on
 $\Delta_2\,:=\,\{\,(x,y)\;|\;x,\,y\,\geq\,0\textrm{\ and\ }x+y\,\leq\,1\;\}$
 \begin{equation}\label{simplex}
  \frac{(\,m+n-2\,)!}{(\,m-2\,)!\,(\,n-2\,)!}
  \;x^{m-2}\;y^{n-2}\;dx\,dy
 \end{equation}
 for all integers $m,\,n\,\geq\,2$ with the characteristic property that
 the moments of the random variable $x^2\mu+y^2\nu$ with arbitrary constants
 $\mu,\,\nu\,\in\,\R$ reproduce the formula of Lemma \ref{mp}
 \begin{eqnarray*}
  \lefteqn{\frac{(\,m+n-2\,)!}{(\,m-2\,)!\,(\,n-2\,)!}\int_{\Delta_2}
  (\,x^2\,\mu\,+\,y^2\,\nu\,)^k\,x^{m-2}\,y^{n-2}\,dx\,dy}\qquad
  &&
  \\
  &=&
  \sum_{r\,=\,0}^k{k\choose r}\,\frac{(\,m+n-2\,)!}{(\,m-2\,)!\,(\,n-2\,)!}
  \frac{(\,m+2r-2\,)!\,(\,n+2(k-r)-2\,)!}{(\,m+n+2k-2\,)!}\;\mu^r\,\nu^{k-r}
 \end{eqnarray*}
 for all $k\,\in\,\N_0$. The proper reason for this striking appearance of
 the simplex $\Delta_2$ in the product formula of Lemma \ref{mp} is of
 course that $\Delta_2$ parametrizes the orbits of the natural action of
 $\mathbf{O}(\,V\,)\,\times\,\mathbf{O}(\,W\,)$ on the Grassmannian
 $\Gr_2(\,V\oplus W\,)$ of $2$--planes in $V\oplus W$. Without discussing
 this rather interesting observation in more detail we use the simplicity
 of the relevant probability measure on $\Delta_2$ to obtain the following
 corollary of Lemma \ref{mp}:

 \begin{Corollary}[Sectional Curvature Densities for Products]
 \hfill\label{scx}\break
  For all $m,\,n\,\geq\,2$ and all $\mu,\,\nu\,\in\,\R$ there exists
  a unique probability measure $\rho^{m,n}_{\mu,\nu}(s)\,ds$ on $\R$
  with support in the convex hull $[\,\min\{\mu,\nu,0\},\,\max\{\mu,\nu,
  0\}\,]$ of $\mu,\,\nu,\,0$ such that the sectional curvature density
  of the product of two compact Riemannian manifolds $M$ and $N$ of
  dimen\-sions $m$, $n$ and sectional curvature densities $\rho^M(\mu)\,d\mu$
  and $\rho^N(\nu)\,d\nu$ can be calculated via:
  $$
   \rho^{M\,\times\,N}(\;s\;)
   \;\;=\;\;
   \int_{\R\,\times\,\R}\rho^{m,\,n}_{\mu,\,\nu}(\,s\,)
   \;\rho^M(\,\mu\,)\;\rho^N(\,\nu\,)\,d\mu\,d\nu
  $$
  The probability density $\rho^{m,n}_{\mu,0}(s)$ for the value $\nu\,=\,0$
  say reads for all $\mu\,\neq\,0$ and $m,\,n\,\geq\,2$
  $$
   \rho^{m,\,n}_{\mu,\,0}(\;s\;)
   \;\;=\;\;
   \frac{(\,m+n-2\,)!}{(\,m-2\,)!\,(\,n-1\,)!}\;\;
   \sqrt{\frac s\mu}^{m-3}\;\Big(\,1\,-\,\sqrt{\frac s\mu}\,\Big)^{n-1}
   \;\frac{ds}{2\,|\,\mu\,|}
  $$
  with support in $[\,0,\,\mu\,]$ or $[\,\mu,\,0\,]$ for positive or negative
  $\mu$ respectively so that $\frac s\mu\,>\,0$ for almost all $s\,\in\,\R$.
  Incidentally this measure can be used for $n\,=\,1$ as well in order to
  calculate the sectional curvature density for the Riemannian product with
  a circle $N\,=\,S^1$.
 \end{Corollary}

 \proof
 For given constants $\mu,\,\nu\,\in\,\R$ and dimensions $m\,,n\,\geq\,2$
 let us consider the function $s\,:=\,x^2\mu+y^2\nu$ on the simplex $\Delta_2$
 as a random variable with respect to the probability measure (\ref{simplex}).
 The inequality $x^2+y^2\,\leq\,x+y\,\leq\,1$ valid for all $(x,y)\,\in\,
 \Delta_2$ implies that $s$ is a convex linear combination of $\mu,\,\nu,\,0$,
 moreover all three values are evidently attained on $\Delta_2$. In consequence
 the random variable $s$ has a probability density $\rho^{m,\,n}_{\mu,\,\nu}(s)
 \,ds$ on the interval $[\,\min\{\mu,\nu,0\},\,\max\{\mu,\nu,0\}\,]$, whose
 moments reproduce the product formula of Lemma \ref{mp}. The moments of the
 probability density $\rho^{M\times N}(s)\,ds$ defined in the Corollary thus
 equal
 \begin{eqnarray*}
  \lefteqn{\int_\R\rho^{M\,\times\,N}(\,s\,)\,s^k\,ds}
  &&
  \\
  &=&
  \int_{\R\,\times\,\R}\Big(\;\int_\R\rho^{m,\,n}_{\mu,\,\nu}(\,s\,)\,
  s^k\,ds\;\Big)\;\rho^M(\,\mu\,)\;\rho^N(\,\nu\,)\;d\mu\,d\nu
  \\
  &=&
  \sum_{r\,=\,0}^k{k\choose r}
  \frac{[\,m+2r-2\,]_{2r}\,[\,n+2(k-r)-2\,]_{2(k-r)}}{[\,m+n+2k-2\,]_{2k}}
  \Big(\int_\R\rho^M(\mu)\,\mu^r\,d\mu\Big)\,
  \Big(\int_\R\rho^N(\nu)\,\nu^{k-r}\,d\nu\Big)
 \end{eqnarray*}
 the moments of the sectional curvature of the Riemannian product $M\times N$
 for all $k\,\in\,\N_0$ so that the corresponding probability measures
 necessarily agree, too.

 In the special case $\nu\,=\,0,\,\mu\,\neq\,0$ the random variable
 $s\,:=\,x^2\mu$ depends on $x$ only so that we may integrate the probability
 density (\ref{simplex}) over $y$ to obtain the distribution of $x$. The
 subsequent substitution $x\,=\,\sqrt{\frac s\mu}$ with $x\,dx\,=\,
 \frac{ds}{2\,|\mu|}$ provides us with the distribution
 $$
  \frac{(m+n-2)!}{(m-2)!\,(n-2)!}
  \,\int_0^{1-x}\!\!\!x^{m-2}\,y^{n-2}\,dy\,dx
  \;\;=\;\;
  \frac{(m+n-2)!}{(m-2)!\,(n-1)!}
  \,\sqrt{\frac s\mu}^{m-3}\,\Big(1-\sqrt{\frac s\mu}\Big)^{n-1}
  \,\frac{ds}{2\,|\,\mu\,|}
 $$
 for the random variable $s$ itself. By assumption $m\,\geq\,2$ so that
 the apparent pole of this density in $s\,=\,0$ is always integrable. A
 discussion of the special case $n\,=\,1$ along a similar line of argument
 is left to the reader.
 \qed

 \pfill
 Let us close this section with a small application illustrating the
 examples and formulas derived above, which exemplifies the use of the
 moments of the sectional curvature in the algebraic framework of local
 Riemannian invariants discussed in \cite{w}. It is well--known that there
 exists a sequence $\pf_1,\,\pf_2,\,\pf_3,\,\ldots$ of natural polynomials
 of degree $1,\,2,\,3,\,\ldots$ respectively for algebraic curvature tensors
 over euclidean vector spaces such that the Euler characteristic of every
 compact Riemannian manifold $M$ of even dimension $m\,=\,2n$ is given by:
 $$
  \chi(\;M\;)\;\;=\;\;\int_M\pf^{T_pM}_n(\;R_p\;)\,\vol_g(\,p\,)
 $$
 In passing we remark that $\pf_n$ is not the better--known Pfaffian
 $\mathrm{Pf}_n$, the latter is not even a natural polynomial in the
 sense discussed in the comments before Corollary \ref{mrm}, nevertheless
 the two polynomials are very closely related. For essentially the same
 reason {\em no} other characteristic number of a compact manifold $M$ can
 be written similarly as the integral of a natural polynomial on curvature
 tensors over $M$, the Pontryagin numbers for example depend explicitly on
 the choice of orientation for $M$.

 Specifically for a homogeneous Riemannian metric $g$ on a manifold $M$
 of even dimension $m\,=\,2n$ the value $\pf^{T_pM}_n(\,R_p\,)\,\in\,\R$
 of this natural polynomial on the curvature tensor $R_p$ is independent
 of the point $p\,\in\,M$ and thus equals the quotient of the Euler
 characteristic $\chi(M)$ of $M$ by its volume with respect to $g$.
 In this way we obtain the explicit values
 $$
  \pf_2(\,R^{S^4}_p\,)
  \;\;=\;\;
  \frac{\chi(\,S^4\,)}{\Vol\;S^4}
  \;\;=\;\;
  \frac{3}{4\,\pi^2}
  \qquad\qquad
  \pf_2(\,R^{\C P^2}_p\,)
  \;\;=\;\;
  \frac{\chi(\,\C P^2\,)}{\Vol\;\C P^2}
  \;\;=\;\;
  \frac{24}{4\,\pi^2}
 $$
 for the round sphere $S^4$ of radius $1$ and volume $\frac{8\pi^2}3$ and
 the complex projective plane $\C P^2$ endowed with the Fubini--Study metric
 $g^\FS$ of volume $\frac1{2\pi}\,\Vol\,S^5\,=\,\frac{\pi^2}2$. A similar
 calculation for the Riemannian product $S^2\times S^2$ of two round spheres
 of radius $1$ and volume $4\pi$ leads to the value $\pf_2(\,R^{S^2\times
 S^2}_p\,)\,=\,\frac1{4\,\pi^2}$. Eventually we consider the Riemannian
 product $S^1\times S^3$ of two round spheres of radius $1$ with vanishing
 Euler characteristic $\pf_2(\,R^{S^1\times S^3}_p\,)\,=\,0$.

 Incidentally three of our four examples, namely $\C P^2,\,S^4$ and $S^2
 \times S^2$ are Einstein manifolds, whereas $S^1\times S^3$ is not. In
 order to discuss this point we extend the scalar product $g$ of a euclidean
 vector space $V$ of dimension $m$ to a scalar product $g^{-1}$ on the space
 $\S^2V^*$ of symmetric $2$--forms on $V$ by choosing an orthonormal basis
 $E_1,\,\ldots,\,E_m$ of $V$ and setting:
 \begin{equation}\label{nq}
  g^{-1}(\;h,\;\tilde h\;)
  \;\;:=\;\;
  \frac12\,\sum_{\mu,\,\nu\,=\,1}^m
  h(\;E_\mu,\;E_\nu\;)\,\tilde h(\;E_\mu,\;E_\nu\;)
 \end{equation}
 This specific normalization of the scalar product $g^{-1}$ is characterized
 by $g^{-1}(\,g,\,g\,)\,=\,\frac m2$ due to the more general $g^{-1}(\,g,\,
 h\,)\,=\,\frac12\,\tr_g\,h\,:=\,\frac12\,\sum h(E_\mu,E_\mu)$ for all
 $h\,\in\,\S^2V^*$, of course other normalizations can and have been
 used in the literature. Regardless of this normalization question the
 scalar product $g^{-1}$ can be used to characterize Einstein manifolds $M$
 of dimension $m$ by the vanishing of the square norm of the so--called
 trace free Ricci tensor
 $$
  |\;\Ric^\circ_p\;|^2
  \;\;:=\;\;
  g^{-1}(\;\Ric_p\,-\,\frac{\kappa_p}m\,g,\;\Ric_p\,-\,\frac{\kappa_p}m\,g\;)
  \;\;=\;\;
  |\;\Ric_p\;|^2\;-\;\frac{\kappa^2_p}{2\,m}
 $$
 where $\Ric_p$ is the Ricci and $\kappa_p\,:=\,\tr_g\,\Ric_p$ the scalar
 curvature of $M$ in $p$. Evidently the norm square of the Ricci tensor is
 additive $|\,\Ric^{M\times N}_{(p,q)}\,|^2\,=\,|\,\Ric_p^M\,|^2\,+\,
 |\,\Ric^N_q\,|^2$  under taking Riemannian products so that the norm
 square of its trace free part $\Ric^\circ_{(p,q)}$ is given by
 $$
  |\;\Ric^\circ_{(p,q)}\;|^2
  \;\;=\;\;
  \frac{\kappa^{M\,2}_p}{2\,m}\;+\;\frac{\kappa^{N\,2}_q}{2\,n}\;-\;
  \frac{(\,\kappa^M_p\,+\,\kappa^N_q\,)^2}{2\,(\,m\,+\,n\,)}
  \;\;=\;\;
  \frac{(\,n\,\kappa^M_p\,-\,m\,\kappa^N_q\,)^2}
  {2\,m\,n\,(\,m\,+\,n\,)}
 $$
 for two Einstein manifolds $M,\,N$ of dimensions $m,\,n\,\geq\,0$, in
 particular $|\,\Ric^\circ\,|\,=\,\frac{(-6)^2}{24}\,\neq\,0$ for the
 Riemannian product $S^1\times S^3$ of two round spheres of radius $1$.
 Completing the values of $\pf_2$ and $|\,\Ric^\circ\,|$ of the four
 examples $\C P^2,\,S^4,\,S^2\times S^2$ and $S^1\times S^3$ by the
 values of their sectional curvature moments calculated in Section \ref{ex}
 and Lemma \ref{mp} respectively we obtain

 \begin{center}
  \begin{tabular}{|c|c|c|c|c|}
   \hline & & & &
   \\[-8pt] 
   \hbox to70pt{\hfil Manifold\hfil}
   & $\quad{\displaystyle\frac43\,\pi^2\,\pf_2}\quad$
   & $\qquad{\displaystyle\Psi_1^2}\qquad$
   & $\qquad{\displaystyle\Psi_2}\qquad$
   & $\qquad{\displaystyle|\,\Ric^\circ\,|}\qquad$
   \\[11pt]
   \hline & & & &
   \\[-6pt]
   $\C P^2$
   & $8$
   & $4$
   & ${\displaystyle\frac{24}5}$
   & $0$
   \\[10pt]
   $S^4$
   & $1$
   & $1$
   & $1$
   & $0$
   \\[10pt]
   $S^2\,\times\,S^2$
   & ${\displaystyle\frac13}$
   & ${\displaystyle\frac19}$
   & ${\displaystyle\frac7{45}}$
   & $0$
   \\[10pt]
   $S^1\,\times\,S^3$
   & $0$
   & ${\displaystyle\frac14}$
   & ${\displaystyle\frac13}$
   & ${\displaystyle\frac32}$
   \\[10pt]
   \hline
  \end{tabular}
 \end{center}
 which essentially proves the following theorem originally due to
 Hitchin--Thorpe \cite{th}:

 \begin{Theorem}[Euler Characteristic of 4--Dimensional Manifolds]
 \hfill\label{ec}\break
  The Euler characteristic of a $4$--dimensional compact Riemannian manifold
  $M$ can be written:
  \begin{eqnarray*}
   \lefteqn{\frac{4\,\pi^2}{3}\;\chi(\,M\,)\;+\;
   \frac49\,\int_M|\;\Ric^\circ_p\;|^2\,\vol_g(\,p\,)}
   \quad
   &&
   \\
   &=&
   \int_M\Psi_2(\,R_p\,)\,\vol_g(\,p\,)\;+\;
   4\,\int_M\Big[\;\Psi_2(\,R_p\,)\;-\;\Psi_1(\,R_p\,)^2\;\Big]\,\vol_g(\,p\,)
  \end{eqnarray*}
  In particular $\chi(M)\,\geq\,0$ for every compact $4$--dimensional
  Einstein manifold $M$ with equality, if and only if $M$ is actually
  a flat Riemannian manifold.
 \end{Theorem}

 \proof
 It is well--known and easy to prove that the vector space of natural
 polynomials of degree $2$ for algebraic curvature tensors over euclidean
 vector spaces has dimension $3$ and is actually spanned by $\kappa^2,\,
 |\,\Ric\,|^2$ and $|\,R\,|^2$. The table of explicit values for the
 four examples $\C P^2,\,S^4,\,S^2\times S^2$ and $S^1\times S^3$ shows
 that the three natural quadratic polynomials $\Psi_1^2,\,\Psi_2$ and
 $|\,\Ric^\circ\,|^2$ are linearly independent and thus span the natural
 quadratic polynomials, in particular we find for the natural quadratic
 polynomial $\frac{4\,\pi^2}3\,\pf_2$ the linear combination
 $$
  \frac{4\,\pi^2}3\,\pf_2
  \;\;=\;\;
  5\,\Psi_2\;-\;4\,\Psi_1^2\;-\;\frac 49\,|\,\Ric^\circ\,|^2
 $$
 again by looking at the explicit values tabulated above. In the resulting
 identity
 $$
  \frac{4\,\pi^2}{3}\;\pf_2(\,R\,)\;+\;\frac49\,|\,\Ric^\circ\,|^2
  \;\;=\;\;
  5\,\Psi_2(\,R\,)\;-\;4\,\Psi_1(\,R\,)^2
  \;\;=\;\;
  \Psi_2(\,R\,)\;+\;4\,\Big[\;\Psi_2(\,R\,)\;-\;\Psi_1(\,R\,)^2\;\Big]
 $$
 valid for all algebraic curvature tensors $R$ the two expressions $\Psi_2(R)$
 and $\Psi_2(R)-\Psi_1(R)^2$ are the second moment and the variance of the
 sectional curvature considered as a random variable on $\Gr_2V$. In particular
 both expressions are non--negative with equality if and only if the random
 variable takes the value $0$ almost everywhere.
 \qed
\end{document}